\documentclass[11pt]{amsart}
\usepackage{amsmath,amsthm, amscd, amssymb, amsfonts}

\def\pf{\begin{proof}}
\def\epf{\end{proof}}

\usepackage{t1enc}
\usepackage{latexsym}
\usepackage{verbatim}
\usepackage{amsgen,amsmath,amstext,amsbsy,amsopn,amsfonts,amssymb}

\newtheorem{maintheorem}{Theorem}
\newtheorem{theorem}{Theorem}[section]
\newtheorem{lemma}[theorem]{Lemma}
\newtheorem{corollary}[theorem]{Corollary}

\newtheorem{proposition}[theorem]{Proposition}

\newtheorem{step}{Step}

\theoremstyle{definition}
\newtheorem{definition}[theorem]{Definition}

\theoremstyle{remark}
\newtheorem{remark}[theorem]{Remark}

\newcommand\ad{\operatorname{ad}}
\newcommand\id{\operatorname{id}}

\newcommand{\rom}{{r}_\Omega}
\newcommand{\rla}{{r}_\Lambda}
\newcommand\re{\operatorname{Re}}
\newcommand{\invo}{{\sigma}}

\newcommand\het{\ell}

\newcommand{\ka}{{\mathfrak k}}
\newcommand{\ko}{{\mathfrak k}_0}
\newcommand{\p}{{\mathfrak p}}
\newcommand{\po}{{\mathfrak p}_0}
\newcommand{\lgot}{{\mathfrak l}}
\newcommand{\m}{\mathfrak{m}}
\newcommand{\ub}{{\mathfrak u}}
\newcommand{\vb}{{\mathfrak v}}
\newcommand{\g}{{\mathfrak g}}
\newcommand{\go}{{\mathfrak g}_0}

\newcommand{\ho}{{\mathfrak h}_0}
\newcommand{\lo}{{\mathfrak l}_0}
\newcommand{\h}{{\mathfrak h}}

\newcommand{\n}{{\mathfrak n}}
\newcommand{\ag}{{\mathfrak a}}
\newcommand{\tg}{{\mathfrak t}}
\newcommand{\ago}{{\mathfrak a}_0}
\newcommand{\agor}{{\mathfrak a}^{\mathbb R}}

\newcommand{\I}{{\tilde J}}
\newcommand{\J}{{\mathbb Z J}}

\newcommand{\dynkin}{\mu}
\newcommand{\tdynkin}{\varsigma_{\mu}}
\newcommand{\tid}{\varsigma}

\newcommand{\jdynkin}{\omega_{\mu, J}}
\newcommand{\jid}{\omega_{J}}

\newcommand{\comp}{\omega}

\newcommand{\ro}{{r}_0}

\newcommand{\Z}{{\mathbb Z}}

\newcommand{\R}{{\mathbb R}}
\newcommand{\Cc}{{\mathbb C}}

\numberwithin{equation}{section}

\newcommand{\End}{\mbox{\rm End\,}}

\begin{document}

\thispagestyle{empty}
\title{A family of Poisson non-compact symmetric spaces}
\author[Andruskiewitsch and Tiraboschi]
{Nicol\'as Andruskiewitsch and Alejandro Tiraboschi}
\address{\noindent
Facultad de Matem\'atica, Astronom\'\i a y F\'\i sica, Universidad
Nacional de C\'ordoba.  CIEM -- CONICET. (5000) Ciudad
Universitaria, C\'ordoba, Argentina}
\email{andrus@famaf.unc.edu.ar, tirabo@famaf.unc.edu.ar}
\thanks{\noindent This work was supported by CONICET, Fund.
Antorchas, Ag. C\'ordoba Ciencia, FONCyT    and Secyt (UNC)}
\date{\today}
\begin{abstract}
We study Poisson symmetric spaces of group type with Cartan
subalgebra ``adapted" to the Lie cobracket.
\end{abstract}
\maketitle

\section*{Introduction}

Let $A$ be a Poisson-Lie group and $T$ a Lie subgroup of $A$. The
homogeneous space $A/T$ endowed with a Poisson structure is a
\emph{Poisson homogeneous space} if the action $A\times A/T \to
A/T$ is a morphism of Poisson manifolds. Poisson homogeneous
spaces, after the seminal paper \cite{Dr1}, have been studied by
several authors, see \cite{EL, FL, L, KRR, K, KoS}; and also in
connection with the quantum dynamical Yang-Baxter equation, see
\cite{EE2, EEM, KS, L2} and references therein.

If $A/T$ carries a Poisson structure such that the natural
projection $A \to A/T$ is a morphism of Poisson manifolds, then
$A/T$ is a Poisson homogeneous space, and in this case is said to
be \emph{of group type}. Assume that $A$  and  $T$ are connected.
Let $\ag$, $\tg$ denote the Lie algebras of $A$, $T$ and let
$\delta: \ag \to \ag \otimes \ag$ be the Lie cobracket inherited
from the Poisson structure of $A$ \cite{D1}, see also \cite[Th.
3.3.1]{KoS}. Then the following conditions are equivalent-- see
\cite{S1, KRR}:
\begin{enumerate}
\item[(i)] $A/T$ is a Poisson homogeneous space of group type;
\item[(ii)] $\{\mu \in \ag^*: \mu(\tg)=0\}$ is a subalgebra in $\ag^*$;
\item[(iii)] $\tg$ is a coideal of $\ag$, i.e. $\delta(\tg) \subset
\tg\otimes \ag + \ag \otimes \tg$.
\end{enumerate}

In this paper, we study Poisson non-compact symmetric spaces of group
type. That is, we assume the following setting:

\begin{itemize}
\item $A = G_0$ is a non-compact
absolutely simple real Lie group with finite center and Lie
algebra $\go$; we fix a Cartan decomposition $\go = \ko \oplus
\po$;
\item the Poisson-Lie group structure on $G_0$ corresponds
to an almost factorizable Lie bialgebra structure on $\go$;
\item $T = K_0$ is a connected Lie subgroup with Lie algebra $\ko$
(in other words, $K_0$ is a maximal compact subgroup of $G_0$).
\end{itemize}

We note that the symmetric space $G_0/K_0$ always has a structure
of Poisson homogeneous space, see subsection \ref{poisshomgt}.
However, whether this Poisson homogeneous structure is of group
type is not evident. Almost factorizable Lie bialgebra structures
on $\go$ were classified in \cite{AJ}, starting from the analogous
classification in the complex case \cite{BD}. In particular, to
each almost factorizable Lie bialgebra structure $\delta$ on $\go$
corresponds a unique Cartan subalgebra $\h$ of the
complexification $\g$ of $\go$, a unique system of simple roots
$\Delta$ in the set of roots $\Phi(\g, \h)$, a unique continuous
parameter $\lambda \in \h^{\otimes 2}$ and a unique
Belavin-Drinfeld triple $(\Gamma_1, \Gamma_2, T)$. Here
$\Gamma_1,$ $\Gamma_2$ are subsets of $\Delta$ and $T: \Gamma_1
\to \Gamma_2$ (see subsection \ref{sec:bialgebras}).

 On the other hand, all maximal compact Lie subgroups of
$G_0$ are conjugated, and they actually arise as the fixed point
set of the Chevalley involution corresponding to some Cartan
subalgebra of $\g$ and some system of simple roots (see subsection
\ref{realforms}). Let $\delta$ an almost factorizable Lie
bialgebra structure  on $\go$, and $\h \subset \g$ the  Cartan
subalgebra and $\Delta$ the set of simple roots determined by
$\delta$. Let $\omega$ be the Chevalley involution that arises
from $\h$ and $\Delta$. We say that $K_0$ is \emph{adapted to
$\delta$} if the Lie algebra of $K_0$ is the fixed point set of
$\omega$.

Now, let $\dynkin: \Delta \to \Delta$ be an automorphism of the
Dynkin diagram of order 1 or 2. Let $J$ be any subset of the set
$\Delta^{\dynkin}$ of simple roots fixed by $\dynkin$. With this
data we can define unique conjugate-linear Lie algebra involutions
$\tdynkin$, $\jdynkin$ of $\g$, see subsection \ref{realforms}
again. It is a well known result that if $\g_0$ is an absolutely
simple real Lie algebra, then it is the set of fixed points of
$\g$, the complexification of $\g_0$, by $\sigma: \g \to \g$  a
conjugate linear involution, where $\sigma$ is $\tdynkin$ or
$\jdynkin$ for some $\mu$ and $J$ (if applies). We denote $\tid =
\tid_{\id}$, $\jid = \comp_{\id, J}$ and $\comp =
\comp_{\id,\Delta}$ (the Chevalley involution). Here is the main
result of the paper.

\begin{maintheorem}\label{main2}
Let $(\go, \delta)$ be an almost factorizable absolutely simple
real Lie bialgebra, let $\sigma$ be the conjugate-linear
involution of $\g$ such that $\g_0 = \g^\sigma$, and let $K_0$ be
the maximal compact Lie subgroup of $G_0$ adapted to $\delta$.
\begin{itemize}
    \item Assume that $\sigma$ is of the form $\tid$, $\tdynkin$ or
$\jid$. Then $G_0/K_0$ is a Poisson homogeneous space of group
type if and only if the Belavin-Drinfeld triple is trivial and
$(\go, \delta)$ is as in Table \ref{tablauno}.

    \item Assume that $\sigma$ is of the form $\jdynkin$ with $\mu\not= \id$ and that
    the Belavin-Drinfeld triple is trivial. Then $G_0/K_0$ is a Poisson homogeneous space of group
type if and only if $\go = \mathfrak{sl}(3,\mathbb R)$ and
$\lambda_{\alpha, \beta} = -\overline{\lambda_{\mu(\alpha),
\mu(\beta)}}$.
\end{itemize}
\end{maintheorem}

\begin{table}[t]
\begin{center}
\begin{tabular}{|l|p{1.2cm}|c|p{3.8cm} |p{4.4cm}|}
\hline

{\bf $\go$} &  $\invo$ & {\bf Type}&   {\bf Continuous\newline
parameter} & {\bf Remarks\vspace{10pt}}

\\\hline
$\g_{\R}$& $\tid$ & all &  $\lambda_{\alpha, \beta} =0$ &
\\ \hline $\mathfrak{su}(n,n+1)$ &    & $A_{2n}$
& &
\\ $\mathfrak{su}(n+1,n+1)$&$\tdynkin,$ & $A_{2n+1}$
&  $\lambda_{\alpha, \beta} =\overline{\lambda_{\mu(\alpha),
\mu(\beta)}}$, &
\\ $\mathfrak{so}(n-1,n+1)$& $\dynkin \neq  \id$ & $D_{n}$
&  $\re(\lambda_{\alpha, \beta} + \lambda_{\alpha, \mu(\beta)}) =
0$ &
\\ $EII$& & $E_{6}$
& &
\\ \hline
 &&&& \textbf{Painted roots:}
\\$\mathfrak{su}(j,n+1-j)$       & $\jid$& $A_n$ & &
$j^{\text{th}}$ root
\\  $\mathfrak{so}(2,2n-1)$ &       & $B_n$  &
$\lambda_{\alpha, \beta} \in i\mathbb{R}$ & first root
\\  $\mathfrak{sp}(n,\mathbb{R})$      &       & $C_n$ &           &  $n^{\text{th}}$ root
\\  $\mathfrak{so}(2,2n-2)$ &       & $D_n$     &          & first root
\\  $\mathfrak{so^*}(2n)$       &       & $D_n$     &       &  $n^{\text{th}}$ root
\\  $EIII$               &       & $E_6$  &        & extreme of the long  branch
\\  $EVII$         &       & $E_7$        & & extreme of  the long branch
 \\
\hline
\end{tabular}
\end{center}

\

\emph{Explanation of the table.} $\invo$ the involution defined by
$\go$, as in \eqref{tdynkin}, \eqref{jdynkin}. The painted roots
are classifiers of Vogan classification. Explanation of Vogan
classification an notations are in \cite{Kn}.

\caption{$G_0/K_0$ Poisson homogeneous space of group type, $K_0$
adapted}\label{tablauno}

\end{table}

Here $\lambda_{\alpha, \beta} \in \mathbb C$ is obtained from
$\lambda - \lambda^{\dag} = \sum_{\alpha, \beta \in \Delta}
\lambda_{\alpha, \beta} h_{\alpha} \wedge h_\beta$, where
$\lambda$ is the continuous parameter and $\lambda^{\dag}$ denotes
the transpose of $\lambda$. The proof of the Theorem follows from
Propositions \ref{invo1} (for $\sigma = \tid$, $\Gamma_1 =
\Gamma_2 = \emptyset$), \ref{invo2} (for $\sigma = \tdynkin$,
$\Gamma_1 = \Gamma_2 = \emptyset$), \ref{invo3} (for $\sigma =
\tid$ or $\tdynkin$, $\Gamma_1 \neq \emptyset$), \ref{invo4} (for
$\sigma = \jid$) and \ref{invo5} (for $\sigma = \jdynkin$,
$\Gamma_1 = \Gamma_2 = \emptyset$), in presence of the information
in \cite[Tables 1.1 and 2.1]{AJ}-- summarized in Proposition
\ref{AJancsa}. The only case that remains open is when $\sigma =
\jdynkin$, $\mu\not= \id$, and the Belavin-Drinfeld triple is
non-trivial.

\medbreak The paper is organized as follows. Section \ref{dos} is
devoted to preliminaries on Lie bialgebras, including the
celebrated theorem of Belavin and Drinfeld, and the classification
result in \cite{AJ}. After this, we prove the main result in
Section \ref{tres}, by a case-by-case analysis.

\section{Lie bialgebras}\label{dos}

\subsection{Simple Lie algebras}\label{simpleliealgebras}

In this section we introduce the notation that will be used
through all the paper. If $\theta$ is a bijection of a set $X$,
then $X^{\theta}$ denotes the fixed-point set of $\theta$. If $a
\in \mathbb C$, we denote by $\overline{a}$ the conjugate of $a$.
We set $i = \sqrt{-1}$. All the Lie algebras in this paper are
finite-dimensional, unless explicitly stated.

 We denote by $\g$ a simple complex Lie algebra and by
$B(\,,\,): \g \times \g \to \mathbb C$ the Killing form on $\g$.
Let $\h$ be a Cartan subalgebra of $\g$. For $\lambda \in h^*$, we
denote $h_\lambda$  the element of $\h$ that satisfies
$B(h_\lambda,h) = \lambda(h)$, for all $h$ in $\h$. We extend
$B(\,,\,)$ to $\h^* \times \h^*$:
$$
B(\lambda,\mu):= B(h_\lambda,h_\mu), \quad \text{for $\lambda,\mu
\in h^*$.}
$$
We denote by $\Phi = \Phi(\g, \h)$ the corresponding root system.
Let $\Delta\subset \Phi$ be a system of simple roots. Let $\Phi^+$
be the set of positive roots with respect to $\Delta$. Given
$\alpha = \sum_{\beta\in \Delta} n_{\beta} \beta\in \Phi^+$, we
denote by $\ell(\alpha) = \sum_{\beta\in \Delta} n_{\beta}$ the
length of $\alpha$.

\medbreak Let $\g_\alpha$ be the root space corresponding to
$\alpha$. Then $\g = \g_+ \oplus \h \oplus \g_-$ is the root space
decomposition of $\g$, where $\g_{\pm} = \oplus_{\alpha \in
\pm\Phi^+} \g_\alpha$.

We choose root vectors $e_\alpha \in \g_\alpha - 0$ such that
\begin{alignat}2
\label{relaciones0} B(e_\alpha, e_{-\alpha})&=1,& &\qquad\text{
for } \alpha \in \Phi.
\end{alignat}
 Then
\begin{alignat}2
\label{relaciones1} [e_\alpha, e_{-\alpha}] &= h_\alpha,  &\qquad
&
\\
\label{relaciones3}[e_\alpha, e_{\beta}] &= 0,& &\text{if }
\alpha+\beta \not= 0 \text{ and }\alpha+\beta \not\in \Phi.
\end{alignat}
For every $\alpha, \beta \in \Phi$ such that $\alpha+\beta
\not=0$, let $N_{\alpha,\beta}$ defined by
\begin{alignat}2
\label{relaciones2}[e_\alpha, e_{\beta}] &=
N_{\alpha,\beta}e_{\alpha+\beta},&\qquad &\text{if } \alpha+\beta
\in\Phi.
\end{alignat}
We set $N_{\alpha,\beta} = 0$, if $\alpha+\beta \not= 0$   and
$\alpha+\beta \not\in \Phi$.
 Thus, we have for all $\alpha, \beta$
and $\gamma$ in $\Phi$ such that $\alpha+\beta \not=0$,
\begin{alignat}2
\label{remarknalphabeta-dos} N_{\alpha,\beta} &= -N_{\beta,\alpha},&&\\
\label{remarknalphabeta-tres} N_{\alpha,\beta} &= N_{\beta,\gamma}
=N_{\gamma,\alpha}, &\qquad& \text{if }\alpha, \beta, \gamma \in
\Phi,\, \alpha+\beta+\gamma=0.
\end{alignat}
The following fact is well-known.
\begin{lemma}
Let $\alpha, \beta \in \Phi$ such that $\alpha + \beta \in \Phi$,
and $\alpha - \beta \not\in \Phi$. Then
\begin{equation}
\label{mualpha5} N_{\alpha,\beta}\, N_{-\alpha,-\beta}=
B(\alpha,\beta)^{-1}.
\end{equation}
\end{lemma}
\begin{proof} If we denote $v = [e_{\alpha},e_{\beta}]$, then
\begin{align*}
[[e_{\alpha},e_{\beta}],[e_{-\alpha},e_{-\beta}]] &=
-[e_{-\alpha},[e_{-\beta},v]]- [e_{-\beta},[v,e_{-\alpha}]],
\\
[e_{-\beta},v] &= -[e_{\alpha}, [e_{\beta},e_{-\beta}]] = - [e_{\alpha}, h_{\beta}]=  B(\alpha,\beta)e_{\alpha} \\
[v,e_{-\alpha}] &=-[[e_{-\alpha},e_{\alpha}],e_{\beta}] =
-[h_{-\alpha},e_{\beta}]= B(\alpha,\beta) e_{\beta}.
\end{align*}
Set $c =\frac1{N_{\alpha,\beta}\, N_{-\alpha,-\beta}}$. Then
\begin{align*}
h_{\alpha+\beta} &=[e_{\alpha+\beta},e_{-\alpha-\beta}] =
\frac1{N_{\alpha,\beta}\, N_{-\alpha,-\beta}}
[[e_{\alpha},e_{\beta}],[e_{-\alpha},e_{-\beta}]]
\\ &= -c\, ([e_{-\alpha},[e_{-\beta},v]]+
[e_{-\beta},[v,e_{-\alpha}]]) =  -c\,B(\alpha,\beta)
([e_{-\alpha},e_{\alpha}]+ [e_{-\beta}, e_{\beta}]) = c\,
B(\alpha,\beta)\, h_{\alpha+\beta},
\end{align*}
and the result follows.
\end{proof}

In the follows, trough all this work, $\g$ will denote a simple
complex Lie algebra and $B(\;,\;)$ the Killing form corresponding
to $\g$. Also, $\h$ denotes a Cartan subalgebra of $\g$, $\Phi =
\Phi(\g, \h)$ the corresponding root system, $\Delta$ a set of
simple roots, $\Phi^+$  the set of positive roots with respect to
$\Delta$, $\g_\alpha$  the root space corresponding to $\alpha \in
\Phi$ and $\g = \g_+ \oplus \h \oplus \g_-$  the root space
decomposition.

\medbreak
\subsection{Absolutely simple real Lie algebras}\label{realforms}

We now describe the real Lie algebras we shall work with. A
finite-dimensional real Lie algebra is called {\it absolutely
simple} if its complexification is a simple complex Lie algebra.
It is well-known that a simple real Lie algebra is either
absolutely simple or is the realification of a complex simple Lie
algebra.

\medbreak Let $\g$ be a simple complex Lie algebra, $\h$ be a
Cartan subalgebra, $\Phi$ be the root system and $\Delta$ be a
system of simple roots. Let $\dynkin: \Delta \to \Delta$ be an
automorphism of the Dynkin diagram of order 1 or 2. We choose
$e_\alpha \in \g_\alpha - 0$ such that $B(e_\alpha,
e_{-\alpha})=1$  for $\alpha \in \Phi$. Let $J$ be any subset of
the set $\Delta^{\dynkin}$ of simple roots fixed by $\dynkin$; let
$\chi_{J}: \Delta \to \{0,1\}$ be the characteristic function of
$J$. Then there exist unique conjugate-linear Lie algebra
involutions $\tdynkin$, $\jdynkin$ of $\g$ given respectively by
\begin{alignat}{2}
\label{tdynkin} \tdynkin (e_{\alpha}) &= e_{\dynkin(\alpha)},
&\qquad \tdynkin (e_{-\alpha}) &= e_{-\dynkin(\alpha)}, \\
\label{jdynkin} \jdynkin (e_{\alpha}) &= (-1)^{\chi_J(\alpha)}
e_{-\dynkin(\alpha)}, &\qquad \jdynkin (e_{-\alpha}) &=
(-1)^{\chi_J(-\alpha)} e_{\dynkin(\alpha)},
\end{alignat}
for all $\alpha \in \Delta$. Necessarily,
\begin{equation}\label{sigmah}\tdynkin (h_{\alpha}) =
h_{\dynkin(\alpha)}, \qquad \jdynkin (h_{\alpha}) =
-h_{\dynkin(\alpha)}.
\end{equation}
We shall write
$$
\tid = \tid_{\id}, \qquad \jid = \comp_{\id, J}, \qquad \comp =
\comp_{\id,\Delta}.
$$
Thus $\comp$ is the Chevalley involution of $\g$, with respect to
$\h$ and $\Delta$,  and the fixed point set of $\comp$ is a
compact form $\ub_0$ of $\g$. Since $\comp(e_\alpha)=
-e_{-\alpha}$ for all $\alpha \in \Phi$-- see \eqref{jdynkingen}
below-- one has:
\begin{equation}\label{compact}
\ub_0 = \sum_{\alpha \in \Phi} \R(i\,h_\alpha) + \sum_{\alpha \in
\Phi} \R(e_\alpha - e_{-\alpha}) + \sum_{\alpha \in \Phi} \R
i\,(e_\alpha + e_{-\alpha}).
\end{equation}

The following lemma -- a variation of \cite[Lemma 2.1]{AJ}-- will
be useful later.

\begin{lemma}\label{involuciones}
Let $\invo: \g \to \g$ be a conjugate-linear Lie algebra
involution such that $\invo (\h) = \h$. Thus, we can define
$\invo^* : \h^*\to \h^*$, the adjoint of $\invo$, by
$\invo^*(\lambda)(h)= \lambda(\invo(h))$, for $\lambda \in H^*$,
$h \in \h$. Then $\invo^* (\Phi) = \Phi$ and $\invo(\g_{\alpha}) =
\g_{\invo^* (\alpha)}$, see \emph{loc. cit}. Then, there exists  a
choice of non-zero root vectors $e_{\alpha}\in \g_{\alpha}$,
$\alpha \in \Phi$,  satisfying \eqref{relaciones0}, such that

\medbreak (a). If $\invo^* (\Delta) = \Delta$, there exists a
unique automorphism $\dynkin: \Delta \to \Delta$ of the Dynkin
diagram of order 1 or 2 (which does not depend on the choice of
the $e_{\alpha}$'s), such that $\invo = \tdynkin$ and
\begin{align}
\label{tdynkingen} \tdynkin (e_{\alpha}) &= e_{\dynkin(\alpha)}, \\
\label{sigmahgen}\tdynkin (h_{\alpha}) &=
h_{\dynkin(\alpha)},\end{align} for all $\alpha \in \Phi$, where
$\dynkin: \Phi \to \Phi$ is the linear extension of $\dynkin$.

\medbreak (b). If $\invo^* (\Delta) = -\Delta$, the there exists a
unique automorphism $\dynkin: \Delta \to \Delta$ of the Dynkin
diagram of order 1 or 2  and a unique subset $\I$ of
$\Phi^{\dynkin}$ (neither $\mu$ nor $\I$ depend on the choice of
the $e_{\alpha}$'s) such that: $\invo = \jdynkin$ with $J = \I
\cap \Delta$. Furthermore,
\begin{align}
\label{jdynkingen} \jdynkin (e_{\alpha}) &= (-1)^{\chi_\I(\alpha)}
e_{-\dynkin(\alpha)}, \\
\label{sigmahgen2}\jdynkin (h_{\alpha}) &= -h_{\dynkin(\alpha)},
\end{align}
for all $\alpha \in \Phi$, where $\dynkin: \Phi \to \Phi$ is the
linear extension of $\dynkin$ and $\chi_\I:\Phi \to \{0,1\}$ the
characteristic function of $\I$.

\medbreak (c). Assume the situation in (b). Let $\alpha, \beta \in
\Phi$ and $N_{\alpha,\beta}$ be defined as in subsection
\ref{simpleliealgebras}. We have
\begin{enumerate}
\item
\begin{align}
\label{chimu} \chi_\I(\alpha) = \chi_\I(-\alpha) &=
\chi_\I(\mu(\alpha)),
\\ \label{mualpha}
(-1)^{\chi_\I(\alpha)+\chi_\I(\beta)}N_{-\mu(\alpha), -\mu(\beta)}
&= (-1)^{\chi_\I(\alpha +\beta)} \overline{N_{\alpha,\beta}}.
\end{align}
\item  Let $\chi_\J: \Phi \to \Z$ be the linear extension of  $\chi_J:\Delta \to \{0,1\}$.
%If $\dynkin = \id$ then
Then
\begin{align}
\label{xik}  (-1)^{\chi_\I(\alpha)} &= (-1)^{\chi_\J(\alpha) +
\ell(\alpha) + 1}
\end{align}
for all $\alpha \in \Phi$.
\end{enumerate}
\end{lemma}

\pf Assume that $\invo^* (\Delta) = \pm \Delta$. Let $\dynkin:
\Delta \to \Delta$ be given by $\dynkin = \pm \invo^*$, according
to the case. Then $\dynkin$ is an automorphism of the Dynkin
diagram, and clearly it has order 1 or 2. Let $f_{\alpha}\in
\g_{\alpha}$, $\alpha \in \Phi$, be any choice of non-zero root
vectors satisfying $B(f_\alpha,f_{-\alpha})=1$. Let $c_{\alpha}
\in \Cc -0$ be such that $\invo(f_{\alpha}) = c_{\alpha}
f_{\invo^*(\alpha)}$, $\alpha \in \Phi$. It is known that $B(x,y)
= \overline{ B(\invo(x) , \invo(y))}$ for all $x$, $y\in \g$, see
\cite[p. 180]{H}. Then, for all $\alpha \in \Phi$, we have $1 =
B(f_{\alpha}, f_{-\alpha}) = \overline{B(\invo(f_{\alpha}) ,
\invo(f_{-\alpha}))} = \overline{c_{\alpha}
c_{-\alpha}B(f_{\invo^*(\alpha)},f_{-\invo^*(\alpha)})}=
\overline{c_{\alpha} c_{-\alpha}}$. Hence, if $\alpha \in \Phi$,
we have
\begin{align}\label{normuno} c_{\alpha} c_{-\alpha} &= 1.
\end{align}
Now, $f_\alpha = \invo^2(f_\alpha) = \invo(c_\alpha
f_{\invo^*(\alpha)}) =
\overline{c_\alpha}\invo(f_{\invo^*(\alpha)}) =
\overline{c_\alpha} c_{\invo^*\alpha}f_\alpha$, thus
\begin{align}
\label{normdos} c_{\alpha} \overline{c_{\invo^*\alpha}} &= 1.
\end{align}

\medbreak We prove (a). In this case $\invo^*=\mu$.  For $\alpha
\in \Phi^+$, define $d_\alpha$, by $c_{\alpha} = \overline{ (
d_{\alpha})^{-1}} d_{\dynkin(\alpha)}$.  The existence of such
$d_{\alpha}$ is clear. For $\alpha \in \Phi$ let $e_{\alpha} =
d_{\alpha}f_{\alpha}$, $e_{-\alpha} =
(d_{\alpha})^{-1}f_{-\alpha}$. Then $B(e_\alpha,e_{-\alpha})=1$
and  $$ \invo (e_{\alpha}) = \overline{d_{\alpha}}c_{\alpha}
(d_{\dynkin(\alpha)})^{-1} e_{\dynkin(\alpha)}=
e_{\dynkin(\alpha)}, \qquad \invo (e_{-\alpha})  =
\overline{(d_{\alpha})^{-1}}c_{-\alpha} (d_{\dynkin(\alpha)})
e_{-\dynkin(\alpha)}= e_{-\dynkin(\alpha)}.
$$ The second formula follows from \eqref{normuno}. The uniqueness of $\dynkin$ is evident and (a) follows.

\medbreak Let us now prove (b). In this case $\invo^*=-\mu$. Let
$\alpha \in \Phi^+$. If $e_{\alpha} = d_{\alpha}f_{\alpha}$,
$e_{-\alpha} = (d_{\alpha})^{-1}f_{-\alpha}$, for  non-zero
scalars $d_{\alpha}$, we have
\begin{equation}\label{eqealfas}
\invo (e_{\alpha}) = \overline{d_{\alpha}}c_{\alpha}
d_{\dynkin(\alpha)} e_{-\dynkin(\alpha)} \quad\text{ and } \quad
\invo (e_{-\alpha}) = \overline{(d_{\alpha})^{-1}}c_{-\alpha}
(d_{\dynkin(\alpha)})^{-1} e_{\dynkin(\alpha)}.
\end{equation}
 \medbreak Assume that $\dynkin(\alpha) \neq \alpha$ with
$\alpha$ positive. Choose $d_{\alpha}$, $d_{\dynkin(\alpha)}$ such
that $c_{\alpha} = (\overline{d_{\alpha}}
d_{\dynkin(\alpha)})^{-1}$. Then $\invo (e_{\alpha}) =
e_{-\dynkin(\alpha)}$, but also $\invo (e_{-\alpha}) =
e_{\dynkin(\alpha)}$ by \eqref{normuno}. Hence $\invo (e_{\pm
\dynkin(\alpha)}) = e_{\mp \alpha}$.

\medbreak  Assume that $\dynkin(\alpha) = \alpha$. In this case we
have $c_{\alpha} \in \R$ by \eqref{normuno} and \eqref{normdos}.
Let $\I = \{\beta \in \Phi: c_{\beta} < 0\}$. For equation
(\ref{eqealfas}) we get
$$
\invo (e_{\pm\alpha}) = |d_{\alpha}|^{\pm 2}c_{\pm\alpha}
e_{\mp\alpha}
$$
Thus, is clear that we can choose $d_{\alpha} \in \R$ such that
$c_{\alpha} d_{\alpha}^2= (-1)^{\chi_\I(\alpha)}$ and it is easy
to see that $\invo (e_{\pm\alpha}) = (-1)^{\chi_\I(\pm\alpha)}
e_{\mp\alpha}$.

Now we see the uniqueness of $\I$: let $f'_\alpha = c f_\alpha$,
then $f'_{-\alpha} = (1/c) f_{-\alpha}$ because the pair
$f'_{\alpha},f'_{-\alpha}$ must satisfy
$B(f'_{\alpha},f'_{-\alpha})=1$. As $\mu(\alpha) = \alpha$, we
have that $\sigma(f'_\alpha) = c'_\alpha f'_{-\alpha}$. On the
other hand $\sigma(f'_\alpha) = \overline{c}\,\sigma(f_\alpha) =
\overline{c}\,c_\alpha\, f_{-\alpha} = \overline{c}\,c \,c_\alpha
\, f'_{-\alpha} = |c|^2\, c_\alpha \, f'_{-\alpha}$. Thus, the
sign of $c_\alpha$ is equal to to the sign of $c'_\alpha$, which
implies that  $\I$ is unique, that is, it does not depend on the
$f_{\alpha}$'s. This shows (b).

\medbreak  We prove (c)(1). The equality \eqref{chimu} is evident;
we obtain \eqref{mualpha} applying $\sigma$ to both sides of the
equation $[e_{\alpha}, e_{ \beta}] = N_{\alpha, \beta
}e_{\alpha+\beta}$.

\medbreak  We prove (c)(2). At the beginning of the proof of this
lemma, we can choose the non-zero root vectors $f_{\alpha}\in
\g_{\alpha}$, $\alpha \in \Phi$ such that if $[f_\alpha, f_\beta]
= M_{\alpha,\beta}f_{\alpha+\beta}$ for all $\alpha, \beta,
\alpha+\beta \in \Phi$,   then $M_{\alpha, \beta}$ are real and
$M_{\alpha, \beta} = -M_{-\alpha, -\beta}$, \emph{cf.} \cite[Th.
6.6]{Kn}.  Let $\gamma, \beta \in \Phi^+$, such that
$\alpha=\gamma+\beta \in \Phi^+$, thus $[e_\gamma, e_\beta] =
N_{\gamma,\beta}e_{\alpha}$ with
$$
N_{\gamma, \beta} = \dfrac{d_\gamma d_\beta}{d_{\alpha}}M_{\gamma,
\beta} \quad \text{ and} \quad N_{-\gamma, -\beta} =
\dfrac{(d_\gamma)^{-1} (d_\beta)^{-1}}{(d_{\alpha
})^{-1}}M_{-\gamma, -\beta},
$$
so $\dfrac{N_{-\gamma, -\beta}}{\overline{N_{\gamma, \beta}}}=C
\dfrac{M_{-\gamma, -\beta}}{M_{\gamma, \beta}}$ with
$C=\dfrac{|d_\alpha|^{2}}{ |d_\gamma|^{2}|d_{\beta}|^{2}}$, where
$|\;|$ denotes the complex norm, so $C>0$. Thus, we can assume
that we have a choice of root vectors as in (b) with $N_{\gamma,
\beta}$ such that
\begin{equation}
\dfrac{N_{-\gamma, -\beta}}{\overline{N_{\gamma, \beta}}} < 0,
\end{equation}
for $\gamma, \beta \in \Phi^+$ with $\gamma+\beta \in \Phi^+$.

We argue then by induction on the length $\het( \alpha)$ of
$\alpha$. Equation \eqref{xik} is evident if $\alpha \in \Delta$.
If $\alpha \in \Phi^+ - \Delta$, fix a decomposition $\alpha =
\gamma + \beta$, where $\gamma\in \Phi^+$ and $\beta \in \Delta$.
Then
\begin{align*}
(-1)^{\chi_\I(\alpha)}e_{-\alpha} &= \invo(e_{\alpha}) =
\dfrac{1}{\overline{N_{\gamma, \beta}}} [\invo(e_{\gamma}),
\invo(e_{\beta})]
\\&= \dfrac{1}{\overline{N_{\gamma, \beta}}} (-1)^{\chi_\J(\gamma) + \het(\gamma) +
\chi_\J(\beta) + \het(\beta)} [e_{-\gamma}, e_{-\beta}] \\&=
\dfrac{N_{-\gamma, -\beta}}{\overline{N_{\gamma, \beta}}}
(-1)^{\chi_\J(\alpha) + \het(\alpha)}e_{-\alpha}.
\end{align*}

Then $\dfrac{N_{-\gamma, -\beta}}{\overline{N_{\gamma, \beta}}} =
- 1$ and (c)(2) follows.

\epf

From completeness we state this theorem of E. Cartan:

\begin{theorem} Let $\go$ be a absolutely simple Lie algebra and $\g$
be the complexification of $\go$. Let $\invo$ be the
conjugate-linear Lie algebra involution such that  $\go
=\g^\invo$. Then, there exists $\h$ a $\sigma$-invariant Cartan
subalgebra of $\g$, a system of simple roots $\Delta$ with
$\sigma^*(\Delta) = \pm \Delta$ and an election of $e_\alpha \in
\g_\alpha-\{0\}$ ($\alpha \in \Phi$); such that $\invo$ is an
involution of the form $\tid, \tdynkin, \jid$ or $\jdynkin$ with
$\dynkin \neq \id$ and $J \subset \Delta^\mu$, but $\invo \neq
\comp$ (if non-compact), or $\go =\g^\comp=\ub_0$ (if compact).
\end{theorem}

\vskip .3cm

Let's denote $\invo$ an involution of the form $\tid, \tdynkin,
\jid$ or $\jdynkin$ with $\dynkin \neq \id$ and $J \subset
\Delta^\mu$, but $\invo \neq \comp$. We denote $\go=\g^\invo$ and
$\ho = \h^\invo$ the set of fixed points by $\invo$ on $\g$ and
$\h$ respectively.

An easy calculation shows that $\invo$ commutes with $\omega$, so
$\omega$ preserves $\go$.   Let $\theta_0: \go \to \go$ be the
linear Lie algebra involution given by the restriction of
$\omega$. Then $\theta_0$ is a Cartan involution of $\go$ and
$\theta_0(\ho)= \ho$. Let $\theta: \g \to \g$ be the
complexification of $\theta_0$; clearly $\theta = \invo\omega$.
The transpose of $\theta$ respect to the Killing form preserves
$\Delta$, and in fact coincides with $\dynkin$. Let
\begin{equation}\label{cartan}
\go = \ko \oplus \po,\quad \text{resp.}\quad \g =\ka \oplus \p,
\end{equation}
be the Cartan decomposition associated to $\theta_0$, resp.
$\theta$.  It is known that $\ko = \ub_0 \cap \go$.  Let $G_0$
(resp. $K_0$) be the connected, simple connected Lie group (resp.
connected Lie subgroup of $G_0$) with Lie algebra $\go$ (resp.
$\ko$).

\medbreak
\subsection{Description of $\ka$} \label{descripciondek}To prove the main result  we will need an explicit expression of $\ka$. In
order to do this, we first write down an explicit form of $\sigma$
and $\theta$.

In the hypothesis of the lemma \ref{involuciones}, if $\invo =
\jdynkin$, we denote $\I$, the unique subset of $\Phi^\mu$
determined by lemma \ref{involuciones}(b). In the follows, we say
that {\it $\I$ is determined by $\jdynkin$}.

Let $\gamma \in \Phi$.
As $\theta= \comp\sigma$ and $\comp(e_\alpha)=-e_{-\alpha}$, we
have, \emph{cf.} \eqref{sigmah}:
\begin{flalign} & \text{If } \sigma= \tdynkin
\text{ then }\theta(e_\gamma)=  - e_{-\mu(\gamma)}, \quad \,
\theta(h_\gamma)=  -h_{\mu(\gamma)}.&\\ \label{dosveinte}
&\text{If } \sigma= \jdynkin \text{ then }\theta(e_\gamma)=
-(-1)^{\chi_\I(\gamma)} e_{\mu(\gamma)}, \quad \theta(h_\gamma)=
h_{\mu(\gamma)}.&
\end{flalign}

\begin{comment}
Recall that $\ko = \ub_0 \cap \go$, where $\ub_0$ is the compact
form of $\g$. So, from formula (\ref{compact}) we can deduce
\begin{equation}\label{k0}
\ko=\sum_{\alpha \in \Phi} \R i (h_\alpha - \sigma(h_\alpha)) +
\sum_{\alpha \in \Phi} \R(e_\alpha - e_{-\alpha}+ \sigma(e_\alpha
- e_{-\alpha})) + \sum_{\alpha \in \Phi} \R i(e_\alpha +
e_{-\alpha}- \sigma(e_\alpha + e_{-\alpha})).
\end{equation}
\end{comment}
Since $\ka = \g^\theta$, an easy calculations shows:

\begin{lemma}\label{generadores}
Let $\sigma$ be an involution of the form
$\tid, \tdynkin, \jid$, or $\jdynkin$ with $\dynkin \neq \id$
and $J \subset \Delta^\mu$. Let $\ka$ and $\ko$ be as in \eqref{cartan}. Then

\begin{flalign} \label{kacaso1}
& \text{If } \sigma= \tid \text{ then }\ka =\langle e_\alpha  -
e_{- \alpha}: \alpha \in \Phi\rangle.&
\\\label{kacaso1bis}
& \text{If } \sigma= \tdynkin \text{ then }\ka =\langle h_\alpha -
h_{\mu(\alpha)}, \,\, e_\alpha  - e_{-\mu(\alpha)} : \alpha \in
\Phi\rangle.&
\\ \label{kacaso2}
&\text{If } \sigma= \jid, \text{let $\I$ the subset of $\Phi^\mu$
determined by $ \jid $.  Then }\ka =\langle h_\alpha , \,\, \big(1
- (-1)^{\chi_\I(\alpha)}\big) e_{\alpha} : \alpha \in
\Phi\rangle.&
\\
\label{kacaso2bis} &\text{If } \sigma= \jdynkin, \text{let $\I$
the subset of $\Phi^\mu$ determined by $\jdynkin$.  Then }\ka
=\langle h_\alpha + h_{\mu(\alpha)}, \,\, e_\alpha -
(-1)^{\chi_\I(\alpha)} e_{\mu(\alpha)} : \alpha \in \Phi\rangle.&
\end{flalign}
\qed \end{lemma}

\medbreak
\subsection{Lie bialgebras}\label{sec:bialgebras}
As is well-known, there is a one-to-one correspondence between
Poisson-Lie group structures on a connected and simply connected
Lie group $A$ and Lie bialgebra structures on its Lie algebra
$\ag$ \cite{D1}. Recall that a (real or complex) Lie bialgebra is
a pair $(\ag, \delta)$ where $\ag$ is a (real or complex) Lie
algebra and $\delta: \ag \to \Lambda^2(\ag)$ is a 1-cocycle
satisfying the co-Jacobi identity.

\medbreak  The Drinfeld double of a Lie bialgebra $\ag$ is denoted
$D(\ag)$-- see \cite{D}. Let $(\ag,\delta)$ be a quasitriangular
Lie bialgebra, \emph{i. e.} there exists $r = \sum_i r_i \otimes
r^i \in \ag \otimes \ag$, such that $\delta = \partial r$ and $r$
satisfies the classical Yang-Baxter equation (CYBE):
$$[r^{12},r^{13}]+[r^{12},r^{23}]+[r^{13},r^{23}]= 0,$$
where if  $r^{12} = \sum_i r_i \otimes r^i \otimes 1$,  $r^{13} =
\sum_i r_i \otimes 1 \otimes r^i $, etc.  In this case, we denote
$(\ag,r)$ instead of $(\ag, \delta)$. We denote $r^{\dag}$ de
transposition of $r$, i.e. $r^\dag =\sum_i  r^i \otimes r_i$. Let
$r_s = r + r^\dag \in S^2 \ag$, then we can define a bilinear form
$\langle\;,\;\rangle_*$ on $\ag^*$ and a map $r_s:\ag^* \to \ag$
by
$$
\langle\mu,\gamma\rangle_* = \sum_i \mu(s_i)\gamma(s^i), \qquad
r_s(\mu) = \sum_i  \mu(s_i)s^i,
$$
where $r_s = \sum_i  s_i \otimes s^i$, $\mu,\gamma \in \ag^*$. We
say that $(\ag,r)$ is {\it factorizable} if
$\langle\;,\;\rangle_*$ is a \emph{nondegenerate} inner product on
$\ag^*$ \cite{RS}.  In this case, the map $r_s$ is an isomorphism
and the bilinear from on $\ag$ defined
$$
\langle x,y \rangle = \langle
r_s(x)^{-1},r_s(y)^{-1}\gamma\rangle_*, \qquad x,y \in \ag
$$
is nondegenerate.  An easy calculation shows that

\begin{proposition}(\cite{RS})\label{amasa}
If $(\ag,r)$ is  a {factorizable} Lie bialgebra, then $D(\ag)$ is
isomorphic to $\ag\oplus\ag$, with bilinear form given by
$$
\langle (x,x'), (y,y')\rangle = \langle x,y \rangle - \langle
x',y'\rangle.
$$
\end{proposition}

\begin{comment}
The \emph{canonical element} of $(\ag,r)$ is the image of $r$
under the Lie bracket: $H_c := \sum_i[r_i, r^i]$, if $r = \sum_i
r_i \otimes r^i$.
\end{comment}

\medbreak  Turning to the real case, a real Lie bialgebra
$(\ago,\delta)$ is {\it almost factorizable} if the
complexification $(\ag,\delta)$ is factorizable (\cite{AJ}). The
following class of almost factorizable Lie bialgebras is important
in the classification of real simple Lie bialgebras given in
\cite{AJ}.

\begin{definition}\label{imaginaryfactorizable}
We say that $(\ago, \delta)$ is \emph{imaginary factorizable} if
the complexification $(\ag, \delta)$ is factorizable and $r\in \ag
\otimes \ag$ is given by
\begin{equation}\label{almost-factorizable}
r = \rla + i \rom, \quad \text{where }\rla \in \Lambda^2(\ago),
\quad \rom \in S^2(\ago).
\end{equation}
\end{definition}
In this case, $D(\ago)$ is isomorphic to the realification
$\agor$, with bilinear form given by $(u \vert v)  = 2\re\langle
u, v\rangle = \langle u, v\rangle - \langle\invo(u),
\invo(v)\rangle$ \cite[Prop. 3.1]{AJ}. Here $\invo$ is the
conjugate-linear involution of $\ag$  such that $\ago =
\ag^{\invo}$.

\medbreak  We review now the celebrated classification of Lie
bialgebra structures on complex simple Lie algebras. Let $\g$ be a
complex semisimple Lie algebra, $\h$ a Cartan subalgebra of $\g$
and $\Delta$ a system of simple roots. Recall that a
\emph{Belavin-Drinfeld triple} (BD-triple for short) is a triple
$(\Gamma_1, \Gamma_2, T)$ where $\Gamma_1,$ $\Gamma_2$ are subsets
of $\Delta$ and $T: \Gamma_1 \to \Gamma_2$ is a bijection that
preserves the inner product coming from the Killing form, such
that the \emph{nilpotency condition} holds: for any $\alpha \in
\Gamma_1$ there exists a positive integer $n$ for which $T^n
(\alpha)$ belongs to $\Gamma_2$ but not to $\Gamma_1$.

\medbreak Let $(\Gamma_1, \Gamma_2, T)$  be a Belavin--Drinfeld
triple. We can extend $T: \Z\Gamma_1 \to \Z\Gamma_2$ in the
natural way. Let $L\subset \Delta$, define $\g_L =
\bigoplus_{\alpha\in \Z L} \g_\alpha \oplus \bigoplus_{\alpha\in
L} \mathbb C h_\alpha$. Then we  can define $\hat T: \g_{\Gamma_1}
\to \g_{\Gamma_2}$ a Lie algebra automorphism, such that $\hat
T(\g_\alpha) = \g_{T(\alpha)}$ for all $\alpha \in \Z\Gamma_1$--
see \cite{BD}.

Let ${\widehat\Gamma}_i = \Z\Gamma_i \cap \Phi^+$ for $i=1,2$.
There is an associated partial ordering on $\Phi^+$ given by
$\alpha \prec \beta$ if $\alpha \in \widehat\Gamma_1,$ $\beta \in
\widehat\Gamma_2,$ and $\beta = T^n(\alpha)$ for an integer $n>0$.

We denote by $\Omega \in \g \otimes \g$ the Casimir element of
$\g$; that is, $\Omega = \sum x_i \otimes x^i$ where $(x_i)$,
$(x^i)$ is any pair of dual basis with respect to the Killing form
of $\g$. We denote by $\Omega_{0}$ the component in $\h \otimes
\h$ of $\Omega$; that is, $\Omega_0 = \sum h_i \otimes h^i$ where
$(h_i)$, $(h^i)$ is any pair of dual basis in $\h$ with respect to
the restriction of the Killing form of $\g$ to $\h$.

\medbreak  A {\it continuous parameter} for the BD-triple
$(\Gamma_1, \Gamma_2, T)$ is an element $\lambda \in \h^{\otimes
2}$ such that
\begin{align}
(T (\alpha) \otimes 1) \lambda + (1 \otimes \alpha) \lambda &= 0,
\quad \text{for all } \alpha \in \Gamma_1,
 \\
\lambda + \lambda^{\dag} &= \Omega_0.
\end{align}
Recall that $\lambda^{\dag}$ denotes the transposition of
$\lambda$.  \medbreak
\begin{theorem}[Belavin--Drinfeld, see \cite{BD}]\label{BDtheorem}
Let $(\g, \delta)$ be  a factorizable complex simple Lie
bialgebra. Then there exist a Cartan subalgebra $\h$, a system of
simple roots $\Delta$, a Belavin--Drinfeld triple $(\Gamma_1,
\Gamma_2, T)$,  a continuous parameter $\lambda$ and $t\in \Cc -
0$ such that the $r$-matrix is given by
\begin{equation}\label{BD}
r = \frac t2 \left(\lambda + \sum_{\alpha \in \Phi^+} x_{-\alpha}
\otimes x_\alpha
        + \sum_{\alpha, \beta \in \Phi^+,  \alpha \prec \beta}
                               x_{-\alpha} \wedge x_\beta\right),
\end{equation}
where $x_\alpha\in \g_\alpha$ are  normalized by
\begin{alignat}2
\label{normalizacion} B(x_\alpha,x_{-\alpha}) &= 1,& \qquad
&\text{for all } \alpha \in
\Phi^+ \\
\hat T(x_\alpha) &= x_{T(\alpha)},&&\text{for all } \alpha \in
\Gamma_1.
\end{alignat}
\end{theorem}

Clearly, $r + r^{\dag} = t\Omega$.  Note that the normalization
condition \eqref{normalizacion} is the same as
\eqref{relaciones0}. Thus, given any family $\{e_{\alpha}: \alpha
\in \Phi\}$ satisfying \eqref{relaciones0}, there exists
$C_{\alpha}\in \mathbb C$ such that

$$x_{\alpha} = C_{\alpha}e_{\alpha}, \qquad
C_{\alpha}C_{-\alpha} = 1.$$

\medbreak We next recall some results of \cite{AJ} about the
classification of real simple Lie bialgebras. Let $\dynkin$ be an
automorphism of the Dynkin diagram.

\begin{itemize}
\item A BD-triple  $(\Gamma_1, \Gamma_2, T)$ is \emph{$\dynkin$-stable}
if $\dynkin(\Gamma_1) =\Gamma_1$, $\dynkin(\Gamma_2) =\Gamma_2$, and
$T\dynkin = \dynkin T$.
\item A BD-triple  $(\Gamma_1, \Gamma_2, T)$ is
\emph{$\dynkin$-antistable} if $\dynkin(\Gamma_1) =\Gamma_2$,
$\dynkin(\Gamma_2) =\Gamma_1$, and $T^{-1}\dynkin = \dynkin T$.
\end{itemize}
If $\dynkin = \id$ then all BD-triples are
$\dynkin$-stable, and the only BD-triple $\dynkin$-antistable has
$\Gamma_1 =  \Gamma_2 = \emptyset$.

\medbreak
\begin{proposition}\label{AJancsa}
 Let $(\go, \delta)$ be an absolutely simple real Lie bialgebra.
Let $\g$ be the complexification of $\go$ and let $\invo$ be the
conjugate-linear involution of $\g$ whose fixed-point set is
$\go$. Assume that $(\go, \delta)$ is almost factorizable. Then
there exist:

\begin{itemize}
\item A Cartan subalgebra $\h$ of $\g$.

\medbreak
\item A system of simple roots $\Delta\subset \Phi(\g, \h)$.

\medbreak
\item A Belavin-Drinfeld triple $(\Gamma_1, \Gamma_2, T)$
and a continuous parameter $\lambda \in \h^{\otimes 2}$.
 Write
$$\lambda - \lambda^{\dag} = \sum_{\alpha, \beta \in \Delta}
\lambda_{\alpha, \beta} h_{\alpha} \wedge h_\beta.$$ By convention
$\lambda_{\alpha, \beta}=-\lambda_{\beta,\alpha}$ for all $\alpha,
\beta \in \Delta$.

\medbreak
\item A complex number $c$ with $c^2 \in \R$; set $t = 2ic$.
\end{itemize}

\medbreak All these data verify:

\medbreak (a) $\h$ is stable under $\invo$ (we denote $\ho := \h
\cap \go$).

\medbreak (b) $\invo^*(\Delta)$ is either $\Delta$ or $-\Delta$;
furthermore $\dynkin: = \invo^*: \Delta \to \pm \Delta$ is an
automorphism of the Dynkin diagram. There are two possibilities:

\begin{enumerate}
\item[{\it(i)}] If $\invo^*(\Delta) = \Delta$ then, by Lemma
\ref{involuciones}, there is an appropriate choice of the
$e_{\alpha}\in \g_\alpha$ ($\alpha \in \Phi$) satisfying
\eqref{relaciones0}, such that $\invo$ is either $\tid$ or
$\tdynkin$ with $\mu\neq \id$ and \eqref{tdynkingen} holds.

\medbreak
If $\invo = \tid$, $t\in \R$, then $\lambda_{\alpha, \beta}\in \R$ for
all $\alpha, \beta \in \Delta$ (no restrictions on the BD-triple).

\medbreak
If $\invo = \tdynkin$, then $t\in \R$, $\lambda_{\alpha,
\beta}= \overline{\lambda_{\mu(\alpha),\mu(\beta)}}$, for all $\alpha,
\beta \in \Delta$ and the BD-triple is $\dynkin$-stable.

\item[{\it (ii)}] If $\invo^*(\Delta) = -\Delta$ then, by  Lemma
\ref{involuciones}, there is an appropriate choice of the
$e_{\alpha}\in \g_\alpha$ ($\alpha \in \Phi$) satisfying
\eqref{relaciones0}, such that $\invo$ is either $\jid$,  $\comp$
or $\jdynkin$ with $\mu\neq \id$ and $J\subset \Delta^\mu$, and
\eqref{jdynkingen} holds.

\medbreak
If $\invo = \jid$, or $\invo = \comp$, then $t\in i\R$,
$\lambda_{\alpha, \beta}\in i\R$, for all $\alpha, \beta \in
\Delta$ and the BD-triple has $\Gamma_1 =  \Gamma_2 = \emptyset$.

\medbreak
If $\invo = \jdynkin$, then $t\in i\R$, $\lambda_{\alpha,
\beta}=-\overline{\lambda_{\mu(\alpha),\mu(\beta)}}$, for all $\alpha,
\beta \in \Delta$ and the BD-triple is $\dynkin$-antistable.
\end{enumerate}

\medbreak (c) $\delta = \partial r$ as in Theorem \ref{BDtheorem}.
Furthermore $\delta = \partial \ro$ where $\ro \in \Lambda^2(\go)$
is given by the formula
\begin{equation}\label{r0}
\ro = \frac t2 \left(\lambda - \lambda^{\dag} + \sum_{\alpha \in
\Phi^+} e_{-\alpha} \wedge e_\alpha
        + \sum_{\alpha, \beta \in \Phi^+,  \alpha \prec \beta}
       C_{-\alpha} C_{\beta} \, e_{-\alpha} \,\wedge
       e_\beta\right),
\end{equation}
with an adequate election of   $C_{\alpha}\in \mathbb C$  for
$\alpha \in \Phi$ such that $C_{\alpha}C_{-\alpha} = 1$.

\medbreak (d)
\begin{equation}\label{titaportitar0}
(\theta\otimes\theta)\ro = \frac{\overline{t}}2
\left(\sum_{\alpha, \beta \in \Delta} \overline{\lambda_{\alpha,
\beta}} h_{\alpha} \wedge h_\beta  - \sum_{\alpha \in \Phi^+}
e_{-\alpha} \wedge e_\alpha + \sum_{\alpha, \beta \in \Phi^+,
\alpha \prec \beta} \overline{C_{-\alpha} C_{\beta}}\, e_{\alpha}
\wedge e_{-\beta}\right).
\end{equation}
\end{proposition}

\pf Parts (a) to (c) are \cite[Lemma 3.1 and Lemma 3.4]{AJ}
combined with Lemma \ref{involuciones}. As $r_0$ belongs to
$\Lambda^2(\g_0)$, $(\sigma\otimes\sigma)r_0= r_0$. Since $\theta=
\comp\sigma$, part (d) follows. \epf

\medbreak
\subsection{Poisson homogeneous spaces of group type arising from graphs}\label{poisshomgt}

\

Let $A$ be a connected and simply connected Poisson-Lie group with
Lie bialgebra $(\ag,\delta)$.  Let $T$ be a connected Lie subgroup
of $A$ with Lie algebra $\tg$. Recall that Poisson homogeneous
structures on $A/T$ are classified by Lagrangian subalgebras
$\lgot$ of $D(\ag)= \ag \oplus \ag^*$, the Drinfeld double of
$\ag$, such that $\lgot \cap \ag = \tg$ \cite{Dr1}.

Recall that  the canonical bilinear form of $D(\ag)$ is given by
$\langle x+\mu|x'+\mu' \rangle = \mu'(x) + \mu(x')$ for $x,x'  \in
\ag$, $\mu,\mu' \in \ag^*$. In this subsection, if $\vb$ is a
subspace of $\ag$, then $\vb^\perp$ denotes the orthogonal
subspace with respect to $\langle \;|\; \rangle$, thus $\vb^\perp
\cap \ag^*= \{\mu \in \ag^*: \mu(\vb)=0\}$ is the annihilator of
$\vb$.

\begin{lemma}\label{lagrangiano} If $\vb$ is a subspace of $\ag$, then $\ub = \vb \oplus (\vb^\perp \cap \ag^*)$ is a
Lagrangian subspace of $D(\ag)$.
\end{lemma}
\begin{proof}
Since $\vb \subset \ag$, $\vb^\perp \cap \ag^* \subset \ag^*$, and
$\ag$, $\ag^*$ are isotropic, we have that $\vb$ and $\vb^\perp
\cap \ag^*$ are isotropic. As $\langle \vb|\vb^\perp \cap
\ag^*\rangle \subset \langle \vb|\vb^\perp \rangle = 0$, $\ub$ is
isotropic. It remains to show that $\ub$ is Lagrangian, or
equivalently that $\dim(\ub) = \dim(\ag) =: n$. Because of the non
degeneracy of the bilinear form on $D(\ag)$ we have that
$\dim(\vb) + \dim(\vb^\perp) = \dim(D(\ag)) = 2n$. But $\vb^\perp
= \ag \oplus (\vb^\perp \cap \ag^*)$, thus $ \dim(\vb^\perp) = n +
\dim(\vb^\perp \cap \ag^*)$. Hence $\dim(\vb) + \dim(\vb^\perp
\cap \ag^*) = n$.
\end{proof}

The following result should be well-known; we give a proof for the
sake of completeness.

\begin{proposition}\label{group-type}
$A/T$ is a Poisson homogeneous space of group type if and only if
there exists a Lagrangian subalgebra $\ub$ of $D(\ag)$ such that
\begin{equation}\label{uno}
\ub = \tg \oplus ( \ub\cap \ag^*).
\end{equation}
\end{proposition}
Note that \eqref{uno} implies $\ub \cap \ag = \tg$.

\begin{proof}
($\Rightarrow$)  Let $\ub = \tg \oplus (\tg^\perp \cap\ag^* )
\subset \tg^\perp$. It is clear that \eqref{uno} holds, and $\ub$
is a Lagrangian subspace by the previous Lemma. It remains to
verify that $\ub$ is a Lie subalgebra. Now, by hypothesis $\tg$
and $\tg^\perp \cap\ag^*=\{\mu \in \ag^*: \mu(\tg)=0\}$ are Lie
subalgebras of $D(\ag)$ (see the Introduction). Let $x \in \tg$,
$y \in \tg^\perp \cap\ag^*$ and $z =[x,y]$.  If $w \in \tg$,
$\langle w|z\rangle = \langle [w,x]|y\rangle \in \langle
\tg|\tg^\perp\rangle=0$. If $w \in \tg^\perp \cap\ag^*$, $\langle
w|z\rangle = -\langle [w,y]|x\rangle \in \langle \tg^\perp
\cap\ag^*|\tg\rangle=0$. Thus $\langle z|\ub\rangle = 0$. Since
$\ub$ is Lagrangian, we conclude that $z \in \ub$.

($\Leftarrow$) If $x \in \tg^\perp \cap\ag^*$, then $\langle x|
\tg\rangle = \langle x|\ag^*\rangle=0$, hence $\langle x|
\ub\rangle =0$. Thus $\langle \tg^\perp \cap\ag^*|\ub\rangle =0$,
and $ \tg^\perp \cap\ag^* \subset \ub$ since $\ub$ is Lagrangian.
Hence $ \tg \oplus (\tg^\perp \cap\ag^*) \subset \ub$. By Lemma
\ref{lagrangiano}, $\tg \oplus (\tg^\perp \cap\ag^* )$ is also a
Lagrangian subspace. Then $ \tg \oplus (\tg^\perp \cap\ag^*)= \ub
$, and this implies that $\tg^\perp \cap\ag^* = \ub \cap\ag^*$ is
a Lie subalgebra of $\ag^* $.
\end{proof}

The following construction of Poisson homogeneous spaces was
observed by C. De Concini, and independiently by Karolinsky
\cite{K}. Let $A$ be a connected (real or complex) Poisson-Lie
group with factorizable  Lie bialgebra  $\ag$: recall that the
Drinfeld double is isomorphic to the Lie algebra $\ag\oplus \ag$,
and the invariant form is given by $\langle (x,x'), (y,y')\rangle
= \langle x,y \rangle - \langle x',y'\rangle$ for $(x,x'), (y,y')
\in \ag\oplus \ag$ (Proposition \ref{amasa}). Let $\rho\in
\operatorname{Aut}(\ag)$ preserving $\langle \,, \,\rangle$. Then,
the graph of $\rho$, namely $\ub_{\rho} = \{(x, \rho(x)): x\in
\ag\},$ is a Lagrangian subalgebra of the Drinfeld double and
$A/T$ is a Poisson homogeneous space, where $T$ is the connected
component of the identity of $A^{\rho}$.

This construction can be extended to the imaginary-factorizable
case. Let $A_0$ be a connected real Poisson-Lie group with
imaginary factorizable Lie bialgebra $(\ago, \delta)$. Let $\ag$
be the complexification of $\ago$ and $\sigma$ the
conjugate-linear automorphism of $\ag$ whose fixed point set is
$\ago$. Let $\theta_0\in \operatorname{Aut}(\ago)$ such that
$\theta :=\theta_0 \otimes \id$ preserves $\langle \,,\,\rangle$.
Let $\mu = \theta \invo= \invo\theta$, a conjugate-linear
automorphism of $\ag$. Let $H$ be the connected component of the
identity of $A_0^{\theta_0}$.

\begin{proposition}\label{group-type-almost} Assume that $\theta_0$ is an involution.
Then $A_0/H$ is a Poisson homogeneous space.
\end{proposition}

\pf As $\ago$ is imaginary factorizable, recall that $D(\ago)$ is
isomorphic to the realification $\agor$, with bilinear form given
by $(u \vert v)  = 2\re\langle u, v\rangle = \langle u, v\rangle -
\langle\invo(u), \invo(v)\rangle$ (see definition
\ref{imaginaryfactorizable} and what follows).  We will show that
the real Lie subalgebra $\m := (\agor)^{\mu}$ is Lagrangian, so,
from the Drinfeld's criterion, $A_0/H$ results a Poisson
homogeneous space. If $u,v\in \m$, then $(u \vert v) =
 \langle u, v\rangle - \langle\invo(u),
\invo(v)\rangle =  \langle u, v\rangle - \langle\theta(u), \theta
(v)\rangle = 0$, thus $\m$ is isotropic. Since $\theta \invo =
\invo \theta$, we have $(\agor)^{\mu}\cap \ago = (\agor)^{ \invo
\theta}\cap \ago = \ago^{\theta_0}$. Also, $\m = \ago^{\theta_0}
\oplus i \po$, where $\po$ is the eigenspace of $\theta_0$ of
eigenvalue $-1$. Thus, $\dim \m = \dim \ago^{\theta_0} + \dim \po
= \dim \ago$, since $\theta_0$ is an involution, and $\m$ is
Lagrangian.  \epf

Note that this Poisson homogeneous space is of group type if and
only if $\m = \ago^{\theta_0} \oplus ( \m\cap r_s(\ago^*))$,
because of Proposition \ref{group-type}. For the definition de
$r_s$, see the beginning of subsection \ref{sec:bialgebras} .

In conclusion, the symmetric spaces $G_0/K_0$ always bear a
structure of Poisson homogeneous space, by Propositions
\ref{AJancsa}, \ref{group-type}-- together with De Concini's
remark-- and \ref{group-type-almost}. In this paper we shall
investigate when $G_0/K_0$ bears a structure of Poisson
homogeneous space \emph{of group type}.

\section{Proof of the main result}\label{tres}

In this section we fix $(\go, \delta)$  an almost factorizable
absolutely simple real Lie bialgebra. Let $G_0$ be a non-compact
absolutely simple real Lie group with finite center and Lie
algebra $\go$ and $K_0$ de maximal compact subgroup of $G_0$
adapted to $\delta$ (see the introduction). As usual, $\ko$
denotes the Lie algebra of $K_0$. Let $\sigma$ be the
conjugate-linear involution of $\g$, the complexification of
$\go$, such that $\g_0 = \g^\sigma$. Now and at the end of this
section we will use the notation of subsection \ref{realforms} and
the notation  and results of proposition \ref{AJancsa}. Also, we
use the description of $\ka$, the complexification of $\ko$, given
in subsection \ref{descripciondek}.

As we said in the Introduction, $G_0/K_0$ is a  a Poisson
homogeneous space of group type  of $G_0$ if and only if $\ko$ is
a coideal of $\go$. Our goal is to determine when $\ko$ is a
coideal of $\go$.

\begin{lemma} \label{coideal}$G_0/K_0$ is a Poisson
homogeneous space of group type  if and only if \, $$\ad \ko
\big((\id - \theta)\otimes(\id - \theta)(\ro) \big)= \ad \ka
\big((\id - \theta)\otimes(\id - \theta)(\ro)\big) = 0.$$
\end{lemma}
\begin{proof} Let $r_0= r_1 + r_2$ with $r_1 \in \go \otimes \ko + \ko \otimes \go$
and $r_2  \in \po \otimes \po$.    If $u \in \ko$, then $\delta
(u) = \ad u\, \ro = \ad u\, r_1 + \ad\, u r_2$ and $\ad u\, r_1
\in \in \go \otimes \ko + \ko \otimes \go$ and $\ad u\, r_2 \in
\po \otimes \po$. Hence $\ko$ is a coideal if and only if $\ad u\,
r_2 =0$ for all $u \in \ko$. Now, if $r_1 = \sum x_i \otimes x^i$
with $x_i$ or $x^i$ in $\ko$  and $r_2 = \sum y_i \otimes y^i$
with $y_i$ and $y^i$ in $\po$, because $\ko$ acts as $\id$ on
$\ko$ and as $-\id$ on $\po$, we have $(\id - \theta)\otimes(\id -
\theta)x_i \otimes x^i =0$ for all $i$, and $(\id -
\theta)\otimes(\id - \theta)y_i \otimes y^i = 4y_i \otimes y^i$.
Thus $(\id - \theta)\otimes(\id - \theta)r_0 = 4r_2$ and $\ad u
\,(\id - \theta)\otimes(\id - \theta)r_0 = 4\ad u \,r_2$. Hence
$\ko$ is a coideal if and only if $\ad u\, r_2 =0$ for all $u \in
\ko$ if and only if $\ad u \,(\id - \theta)\otimes(\id -
\theta)r_0 =0$  for all $u \in \ko$ .
\end{proof}

We next give an explicit expression of
$$\widetilde{\ro}:=(\id - \theta)\otimes(\id - \theta)(\ro)
=\ro + (\theta\otimes\theta)(\ro)
 - (\id\otimes \theta + \theta\otimes \id)(\ro)$$
according to the different possibilities for $\invo$.
Then we analyze when $\ko$ is a coideal case by case.

\medbreak
\subsection{Computation of $(\id - \theta)\otimes(\id - \theta)(\ro)$}

\

In the calculations below, keep in mind the equality $(f\otimes
\id + \id \otimes f)(a \wedge b) = f(a) \wedge b+a \wedge f(b)$,
$a, b \in V$, $f\in \End V$. Set
\begin{align}\label{tab}
t_{\alpha,\beta}&=2 Re(\lambda_{\alpha, \beta}
+\lambda_{\alpha,\mu(\beta)}), &&\alpha, \beta \in \Delta,
\\\label{sab}
s_{\alpha,\beta} &=2\,i\,Im(\lambda_{\alpha, \beta}-
\lambda_{\alpha,\mu(\beta)}),  &&\alpha, \beta \in \Delta,
\\\label{dab}
d_{\alpha,\beta} &= C_{-\alpha} C_{\beta}, &&\alpha, \beta \in
\Phi^+, \quad \alpha \prec \beta.\end{align}

\medbreak

\begin{proposition} If $\sigma=\tid$, then
\begin{align}\label{tid-deltau}
\widetilde{\ro} =
 \frac t2\left(
\sum_{\alpha, \beta \in \Delta} 2\lambda_{\alpha,\beta}\,
h_{\alpha} \wedge h_\beta   +\sum_{\alpha \prec \beta}
(\overline{d_{\alpha,\beta}}\,e_{-\alpha} \wedge e_{\beta}
+d_{\alpha,\beta}\, \big(e_{\alpha} \wedge e_{-\beta}+ e_{\alpha}
\wedge e_\beta+ e_{-\alpha} \wedge e_{- \beta}\big))\right).
\end{align}

\medbreak
If $\sigma=\tdynkin$, $\mu \neq \id$, then
\begin{alignat}2 \label{tdynkin-deltau}
\widetilde{\ro} &=&
 \frac t2&\left(
  \sum_{\alpha, \beta \in \Delta} t_{\alpha,\beta}\,
h_{\alpha} \wedge h_\beta  +\sum_{\alpha \prec \beta}
\overline{d_{\alpha,\beta}}\, e_{\alpha} \wedge e_{-\beta} +
d_{\alpha,\beta}\, \big(e_{-\alpha} \wedge e_{\beta} +
e_{\mu(\alpha)} \wedge e_\beta +
 e_{-\alpha} \wedge e_{- \mu(\beta)}\big))\right).
\end{alignat}

\medbreak If $\sigma=\jid$, then
\begin{equation}\label{jid-deltau}
\widetilde{\ro} = t \sum_{\alpha\in \Phi^+} (1 +
(-1)^{\chi_\I(\alpha)}) e_{-\alpha} \wedge e_\alpha.
\end{equation}

\medbreak If $\sigma=\jdynkin$, $\mu \neq \id$, then
\begin{multline} \label{jdynkin-deltau}
\widetilde{\ro}  = \frac t2\Big(\sum_{\alpha, \beta \in \Delta}
s_{\alpha,\beta}\,h_{\alpha}\wedge h_\beta + 2\sum_{\alpha\in
\Phi^+} (e_{-\alpha} \wedge e_\alpha + (-1)^{\chi_\I(\alpha)}
e_{-\mu(\alpha)}\wedge e_\alpha)
\\ + \sum_{\alpha \prec \beta}\Big(-\overline{d_{\alpha,\beta}}\,
e_{\alpha} \wedge e_{-\beta} +d_{\alpha,\beta}\, \big(e_{-\alpha}
\wedge e_{\beta} + (-1)^{\chi_\I(\alpha)} \, e_{-\mu(\alpha)}
\wedge e_\beta + (-1)^{\chi_\I(\beta)} \, e_{-\alpha} \wedge e_{
\mu(\beta)}\big)\Big)\Big).
\end{multline}
\end{proposition}

\pf We probe simultaneously (2.4) and (2.5). For what follows we use the expressions of $\ro$ and
$(\theta\otimes \theta) (\ro)$ given by (\ref{r0}) and
(\ref{titaportitar0}), respectively. Assume that $\invo =
\tdynkin$ with $\mu$ arbitrary. In this case, $t \in \R$, so
$$
\ro + (\theta\otimes\theta)(\ro)  = \frac t2 \left( \sum_{\alpha,
\beta \in \Delta} 2Re(\lambda_{\alpha, \beta}) h_{\alpha} \wedge
h_\beta  + \sum_{\alpha, \beta \in \Phi^+, \alpha \prec \beta}
(d_{\alpha,\beta}\, e_{-\alpha} \wedge e_{\beta} +
\overline{d_{\alpha,\beta}}\, e_{\alpha} \wedge e_{-\beta})\right)
$$
Now,

\begin{alignat*}2
(\theta\otimes \id + \id \otimes\theta) (\ro) &=&
 \frac t2& \Big( \sum_{\alpha, \beta \in \Delta}
 \lambda_{\alpha, \beta} (\theta(h_{\alpha}) \wedge h_\beta+h_{\alpha} \wedge \theta(h_\beta))
 \\
 &&+& \sum_{\alpha\in \Phi^+} (\theta(e_{-\alpha}) \wedge e_\alpha+e_{-\alpha} \wedge
 \theta(e_\alpha))\\
 &&+&\sum_{\alpha, \beta \in \Phi^+,  \alpha \prec \beta}
 d_{\alpha,\beta}\,(\theta(e_{-\alpha}) \wedge e_\beta+e_{-\alpha} \wedge
 \theta(e_\beta)) \Big)\\
 &=&
 -\frac t2& \Big( \sum_{\alpha, \beta \in \Delta}
 \lambda_{\alpha, \beta} (h_{\mu(\alpha)} \wedge h_\beta+h_{\alpha} \wedge h_{\mu(\beta)})
+ \sum_{\alpha\in \Phi^+}  ( e_{\mu(\alpha)}\wedge e_\alpha+
  e_{-\alpha}\wedge e_{-\mu(\alpha)}) \\
 &&+&\sum_{\alpha, \beta \in \Phi^+,  \alpha \prec \beta}
d_{\alpha,\beta}( e_{\mu(\alpha)} \wedge e_\beta + e_{-\alpha}
\wedge
 e_{- \mu(\beta)})\Big).
\end{alignat*}

\medbreak

Here the first term is
$$
\sum_{\alpha, \beta \in \Delta}\lambda_{\alpha, \beta}
(h_{\mu(\alpha)} \wedge h_\beta+h_{\alpha} \wedge h_{\mu(\beta)}) =
\sum_{\alpha, \beta \in \Delta}(\lambda_{\mu(\alpha), \beta} +
\lambda_{\alpha,\mu(\beta)}) h_{\alpha} \wedge h_\beta
= \sum_{\alpha, \beta \in \Delta}2
Re(\lambda_{\alpha,\mu(\beta)}) h_{\alpha} \wedge h_\beta;
$$
the first equality by a change of variables, the second because
in this case
$\lambda_{\alpha,\beta}=\overline{\lambda_{\mu(\alpha),\mu(\beta)}}$,
thus $\lambda_{\mu(\alpha),\beta}=
\overline{\lambda_{\alpha,\mu(\beta)}}$, and thus
$\lambda_{\mu(\alpha), \beta} + \lambda_{\alpha,\mu(\beta)}
= 2 Re(\lambda_{\alpha,\mu \beta})$.

\medbreak On the other hand, the second term is $ \sum_{\alpha\in
\Phi^+}  e_{\mu \alpha}\wedge e_\alpha =\sum_{\alpha\in \Phi^+,
\alpha\not=\mu(\alpha)}   e_{\mu(\alpha)}\wedge e_\alpha$. We can
enumerate the set of those roots such that
$\alpha\not=\mu(\alpha)$ in this way: $\alpha_1, \mu(\alpha_1),
\ldots , \alpha_k, \mu(\alpha_k)$, where  $\alpha_i \not=
\mu(\alpha_i)$, $\alpha_i \not= \alpha_j, \mu(\alpha_j)$, for all
$i\not=j$. Then $ \sum_{\alpha\in \Phi^+}   e_{\mu(\alpha)}\wedge
e_\alpha = \sum_{i}  e_{\mu(\alpha_i)}\wedge e_{\alpha_i}+ e_{
\alpha_i}\wedge e_{\mu(\alpha_i)} = 0$. In analogous way,
$\sum_{\alpha\in \Phi^+} e_{-\alpha}\wedge e_{-\mu(\alpha)}=0$.
So,
\begin{align*}
(\theta\otimes \id + \id \otimes\theta)(\ro) =  -\frac t2 \Big(
\sum_{\alpha, \beta \in \Delta}2 Re(\lambda_{\alpha,\mu(\beta)})
h_{\alpha} \wedge h_\beta + \sum_{\alpha, \beta \in \Phi^+, \alpha
\prec \beta} d_{\alpha,\beta}(e_{\mu(\alpha)} \wedge e_\beta +
e_{-\alpha} \wedge e_{- \mu(\beta)})\Big).
\end{align*}
Hence $\widetilde{\ro} =\ro + (\theta\otimes\theta)(\ro)
 - (\id\otimes\theta + \theta\otimes \id)(\ro)$
equals the right-hand side of \eqref{tdynkin-deltau}; when $\mu =
\id$, this reduces to \eqref{tid-deltau}, by Proposition
\ref{AJancsa} (b) (i).

\medbreak
Now, we probe  (2.6) and (2.7).  Assume that $\invo = \jdynkin$ with $\mu$ arbitrary.
 In this case $t \in i\R$, so
\begin{align*}
\ro + (\theta\otimes\theta)\ro = \frac{t}2 \Big(2i\sum_{\alpha,
\beta \in \Delta} Im{\lambda_{\alpha, \beta}}\, h_{\alpha} \wedge
h_\beta  &+ 2\sum_{\alpha \in \Phi^+} e_{-\alpha} \wedge e_\alpha
\\ &+ \sum_{\alpha, \beta \in \Phi^+, \alpha \prec \beta}
(d_{\alpha,\beta}\,e_{-\alpha} \wedge
e_{\beta}-\overline{d_{\alpha,\beta}}\,e_{\alpha} \wedge
e_{-\beta})\Big).
\end{align*}
Now,
\begin{alignat*}2
(\theta\otimes \id + \id \otimes\theta)\ro &=& \frac t2& \Big(
\sum_{\alpha, \beta \in \Delta}\lambda_{\alpha, \beta}
(\theta(h_{\alpha}) \wedge h_\beta+h_{\alpha} \wedge
\theta(h_\beta))
 +  \sum_{\alpha\in \Phi^+} (\theta(e_{-\alpha}) \wedge e_\alpha+e_{-\alpha} \wedge
 \theta(e_\alpha))\\
 &&+&\sum_{\alpha, \beta \in \Phi^+,  \alpha \prec \beta}
 d_{\alpha,\beta}(\theta(e_{-\alpha}) \wedge e_\beta+e_{-\alpha} \wedge
 \theta(e_\beta))\Big)
\\ &=&
\frac t2& \Big( \sum_{\alpha, \beta \in \Delta}\lambda_{\alpha,
\beta} (h_{\mu(\alpha)} \wedge h_\beta+h_{\alpha} \wedge
h_{\mu(\beta)}) \\
&&-& \sum_{\alpha\in \Phi^+} (-1)^{\chi_\I(\alpha)} (
e_{-\mu(\alpha)}\wedge e_\alpha+
e_{-\alpha}\wedge e_{\mu(\alpha)})\\
&&-& \sum_{\alpha, \beta \in \Phi^+,  \alpha \prec \beta}
d_{\alpha,\beta}\,((-1)^{\chi_\I(\alpha)}\, e_{-\mu(\alpha)}
\wedge e_\beta + (-1)^{\chi_\I(\beta)} \, e_{-\alpha} \wedge e_{
\mu(\beta)})\Big).
\end{alignat*}

 In this case
$\lambda_{\alpha,\beta}=-\overline{\lambda_{\mu(\alpha),\mu(\beta)}}$,
then $\lambda_{\mu(\alpha),\beta}=
-\overline{\lambda_{\alpha,\mu(\beta)}}$, thus $\lambda_{\mu(\alpha),
\beta} + \lambda_{\alpha,\mu(\beta)} = 2\,i\,
Im(\lambda_{\alpha,\mu(\beta)})$ and
$$
\sum_{\alpha, \beta \in \Delta}\lambda_{\alpha, \beta}
(h_{\mu(\alpha)} \wedge h_\beta+h_{\alpha} \wedge h_{\mu(\beta)}) =
\sum_{\alpha, \beta \in \Delta}2\, i\, Im(\lambda_{\alpha,\mu
\beta})h_{\alpha} \wedge h_{\beta}.
$$

Now the second term is
\begin{align*}
\sum_{\alpha\in \Phi^+} (-1)^{\chi_\I(\alpha)} (
e_{-\mu(\alpha)}\wedge e_\alpha+ e_{-\alpha}\wedge
e_{\mu(\alpha)})  = 2\sum_{\alpha\in \Phi^+}
(-1)^{\chi_\I(\alpha)} \,e_{-\mu(\alpha)}\wedge e_\alpha.
\end{align*}
 Thus
\begin{align*}
(\theta\otimes \id + \id \otimes\theta)(\ro) &=
 \frac t2 \Big(2i\sum_{\alpha, \beta \in \Delta}
Im(\lambda_{\alpha,\mu \beta})h_{\alpha} \wedge h_{\beta} -
2\sum_{\alpha\in \Phi^+} (-1)^{\chi_\I(\alpha)}\,
e_{-\mu(\alpha)}\wedge e_\alpha
\\ &\quad -  \sum_{\alpha, \beta \in \Phi^+,  \alpha \prec \beta}
d_{\alpha,\beta}\,((-1)^{\chi_\I(\alpha)}\, e_{-\mu(\alpha)}
\wedge e_\beta+(-1)^{\chi_\I(\beta)}\, e_{-\alpha} \wedge e_{
\mu(\beta)})\Big).
\end{align*}
Hence $\widetilde{\ro} =\ro + (\theta\otimes\theta)(\ro)
 - \big( (\id\otimes\theta)+ (\theta\otimes \id)\big)(\ro)$
equals the right-hand side of \eqref{jdynkin-deltau}. When $\mu =
\id$, this reduces to \eqref{jid-deltau}; indeed recall that
$\Gamma_1 = \Gamma_2 = \emptyset$ in this case by \cite{AJ}. \epf

\medbreak
\subsection{Case $\sigma= \tdynkin$, $\Gamma_1 = \Gamma_2 = \emptyset$}
We begin with arbitrary $\mu$. Recall the definition of $t_{\alpha,\beta}$ in \eqref{tab}.

\begin{proposition} $\ko$ is a coideal of $\go$ if and only if
$t_{\alpha,\beta}=0$ for all $\alpha,\beta \in \Delta$.
\end{proposition}
\begin{proof} Recall that in this case $\ka =\langle h_\alpha -
h_{\mu(\alpha)}, \,\, e_\alpha  - e_{-\mu(\alpha)} : \alpha \in
\Phi\rangle$.  Here $\widetilde{\ro} =
 \dfrac t2  \sum_{\alpha, \beta \in \Delta} t_{\alpha,\beta} h_{\alpha} \wedge h_\beta$ by \eqref{tid-deltau} or \eqref{tdynkin-deltau}.
Note that $t_{\beta,\alpha} =- t_{\alpha,\beta}$. In particular
$\ad \h(\widetilde{\ro})=0$. Thus, we only must calculate
$\ad(e_\gamma -   e_{-\mu(\gamma)})\widetilde{\ro}$, for all
$\gamma \in \Phi$. Now:

\begin{align*}
-\ad e_\gamma \Big(\sum_{\alpha, \beta \in \Delta} t_{\alpha,\beta}
h_{\alpha} \wedge h_\beta \Big)&=
 \sum_{\alpha, \beta \in\Delta}
t_{\alpha,\beta}\big( B(\gamma,\alpha) \, e_\gamma  \wedge h_\beta
+ B(\gamma,\beta)h_{\alpha}\,\wedge e_\gamma\big)\\
&=
 \sum_{\alpha, \beta \in
\Delta} t_{\alpha,\beta} B(\gamma,\alpha) e_\gamma  \wedge h_\beta
+ \sum_{\alpha, \beta \in
\Delta} t_{\beta,\alpha}B(\gamma,\alpha)h_{\beta} \wedge e_\gamma\\
&=
 \sum_{\alpha, \beta \in
\Delta} t_{\alpha,\beta} B(\gamma,\alpha) e_\gamma  \wedge h_\beta
- \sum_{\alpha, \beta \in
\Delta} t_{\alpha,\beta}B(\gamma,\alpha)h_{\beta} \wedge e_\gamma\\
&=
 \sum_{\alpha, \beta \in
\Delta} 2t_{\alpha,\beta} B(\gamma,\alpha) e_\gamma  \wedge
h_\beta =
 e_\gamma  \wedge \sum_{ \beta \in
\Delta}  B(\gamma,\sum_{\alpha \in
\Delta}2t_{\alpha,\beta}\alpha) h_\beta.
\end{align*}
So
$$
\ad(e_\gamma  -   e_{-\mu(\gamma)}) (\widetilde{\ro}) =-\frac t2
e_\gamma \wedge \sum_{ \beta \in \Delta}  B(\gamma,\sum_{\alpha
\in \Delta}2t_{\alpha,\beta}\alpha) h_\beta + \frac t2
e_{-\mu(\gamma)} \wedge \sum_{ \beta \in \Delta}
B(-\mu(\gamma),\sum_{\alpha \in \Delta}2t_{\alpha,\beta}\alpha)
h_\beta.
$$
The terms in the right-hand side are linearly independent, so
$[e_\gamma  -   e_{-\mu(\gamma)}, \widetilde{\ro}]=0$ if and only
if $B(\gamma,\sum_{\alpha \in \Delta}2t_{\alpha,\beta}\alpha)=0$
for all $\gamma\in \Phi^+$, $\beta \in \Delta$. Because of the
non-degeneracy of $B(\,,\,)$, $\ad \ka \, \widetilde{\ro} = 0$ if
and only if $\sum_{\alpha \in \Delta}2t_{\alpha,\beta}\alpha=0$
for all $\beta \in \Delta$, if and only if $t_{\alpha,\beta}=0$
for all $\alpha,\beta \in \Delta$.
\end{proof}

The following propositions follow immediately from the previous
one and \cite[Tables 1.1 and 2.1]{AJ}, see Proposition
\ref{AJancsa}.

\medbreak
\begin{proposition}\label{invo1}If $\invo= \tid$ and
$\Gamma_1 = \Gamma_2 = \emptyset$, then $\ko$ is a coideal of
$\go$ if and only if $\lambda_{\alpha,\beta}=0$ for all
$\alpha,\beta \in \Delta$. \qed\end{proposition}

\medbreak
\begin{proposition}\label{invo2} If $\invo= \tdynkin$ with $\mu \neq \id$ and
$\Gamma_1 = \Gamma_2 = \emptyset$, then $\ko$ is a coideal of
$\go$ if and only if the continuous parameter $\lambda$ satisfies
$\lambda_{\alpha,\beta}=
\overline{\lambda_{\mu(\alpha),\mu(\beta)}}$ and
$t_{\alpha,\beta}=0$, for all $\alpha,\beta \in \Delta$.
\qed\end{proposition}

\medbreak
\subsection{Case $\sigma= \tdynkin$, $\Gamma_1 \not =
\emptyset$,  $\Gamma_2 \not = \emptyset$.}
For arbitrary $\mu$, we have:

\begin{proposition}\label{invo3}  $\ko$ is not a coideal of $\go$.
\end{proposition}

\begin{proof}
Recall that  $\ka$ is generated by $h_{\gamma}-h_{\mu(\gamma)}$,
$e_\gamma - e_{-\mu(\gamma)}$, $\gamma \in \Phi$, see Lemma
\ref{generadores}. We study the action of
$h_{\gamma}-h_{\mu(\gamma)}$, $\gamma \in \Phi$, on
$\widetilde{\ro}$:
\begin{align*}
\ad(h_{\gamma}-h_{\mu(\gamma)}) (\widetilde{\ro})&=\frac t2
\Big(\sum_{\alpha \prec \beta} \,(B(\gamma
-\mu(\gamma),-\alpha+\beta)
\,\overline{d_{\alpha,\beta}}\,e_{-\alpha} \wedge e_{\beta} +
B(\gamma -\mu(\gamma),\alpha-\beta) d_{\alpha,\beta}e_{\alpha}
\wedge e_{-\beta}
\\ + & B(\gamma -\mu(\gamma),\mu(\alpha)+\beta)  \,d_{\alpha,\beta}\, e_{\mu
\alpha} \wedge e_\beta+ B(\gamma -\mu(\gamma),-\alpha-\mu(\beta))
\,d_{\alpha,\beta}\, e_{-\alpha} \wedge
 e_{- \mu(\beta)})\Big)
\\&= \frac t2
\sum_{\alpha \prec \beta}  \,B(\gamma
-\mu(\gamma),-\alpha+\beta)\, A_{\alpha,\beta},
\end{align*}
where $A_{\alpha,\beta}= \overline{d_{\alpha,\beta}} e_{-\alpha}
\wedge e_{\beta} - d_{\alpha,\beta} \,
 e_{\alpha} \wedge e_{-\beta}+   \,d_{\alpha,\beta}\,
e_{\mu(\alpha)} \wedge e_\beta+  \,d_{\alpha,\beta}\, e_{-\alpha}
\wedge
 e_{- \mu(\beta)}$.
Here in the second equality we have argued as follows:
$B(\alpha+\mu(\alpha),\gamma -\mu(\gamma))= 0$, because $\mu^2=
\id$ and $\mu$ is $B$-invariant. Then $B(\gamma
-\mu(\gamma),\mu(\alpha)+\beta)= B(\gamma
-\mu(\gamma),-\alpha+\beta)$. Also, $B(\gamma
-\mu(\gamma),-\alpha-\mu(\beta))= B(\gamma
-\mu(\gamma),-\alpha+\beta)$.

\medbreak It is clear that the elements $A_{\alpha,\beta}$ are
linearly independent, \emph{e. g.} projecting them to $\g^+\otimes
\g^-$ along the root space decomposition. Thus,
$[h_{\gamma}-h_{\mu(\gamma)}, \widetilde{\ro}]= 0$ if and only if
$B(\gamma -\mu(\gamma),\alpha-\beta)=0$, for all $\gamma\in\Phi,
\alpha\in\widehat\Gamma_1, \beta \in \widehat\Gamma_2,
\alpha\prec\beta$. Now, $B(\gamma -\mu(\gamma),\alpha-\beta)=0$ if
and only if $B(\gamma
,\alpha-\beta)=B(\mu(\gamma),\alpha-\beta)=B(\gamma,\mu(\alpha-\beta))$,
for all $\gamma\in \Phi$ if and only if $\alpha-\beta =
\mu(\alpha-\beta)$ if and only if $\alpha +\mu(\beta) =
\mu(\alpha) +\beta$. As $\alpha \prec \beta$, we have that $\alpha
\not=\beta$. If $\alpha\in \Gamma_1$ and $\beta \in \Gamma_2$,
then $\alpha,\mu(\alpha),\beta,\mu(\beta) \in \Delta$ and $\alpha
+\mu(\beta) = \mu(\alpha) +\beta$ if and only if $\alpha=
\mu(\alpha)$ and $\beta= \mu(\beta)$. We conclude that
\begin{equation}
\ad(h_{\gamma}-h_{\mu(\gamma)}) (\widetilde{\ro})= 0,\,\, \forall
\gamma\in \Phi \quad \text{if and only if} \quad \alpha=
\mu(\alpha),\,\, \forall \alpha\in \Gamma_1 \cup\Gamma_2.
\end{equation}

In view of this, we shall assume in the rest of the proof that $\alpha=
\mu(\alpha),\,\, \forall \alpha\in \Gamma_1 \cup\Gamma_2$. Then,
\begin{alignat}2 \label{r0caso1}
\widetilde{\ro} &=&
 \frac t2&
  \sum_{\alpha, \beta \in \Delta} t_{\alpha,\beta} h_{\alpha} \wedge h_\beta
  +\frac t2 \sum_{\alpha \prec \beta} \,\overline{d_{\alpha,\beta}}\,e_{-\alpha} \wedge
e_{\beta} +\,d_{\alpha,\beta}\,e_{\alpha} \wedge e_{-\beta}+
\,d_{\alpha,\beta}\,e_{\alpha} \wedge e_\beta+ \,d_{\alpha,\beta}\,
d _\beta e_{-\alpha} \wedge
 e_{- \beta}.
\end{alignat}

Now, reorder the simple roots in the following way: let $k =
\#(\Gamma_1\cup\Gamma_2)$, and let $\alpha_{k+1},\ldots,\alpha_n$
the simple roots that not belong to $\Gamma_1\cup\Gamma_2$. Let
$\{\beta_1,\ldots, \beta_t\}$ the roots in $\Gamma_1$ that are not
image of $T$, then
\begin{equation} \label{ordenraices}
(\alpha_1, \ldots, \alpha_k) =(\beta_1, T(\beta_1),
\ldots,T^{s_1}(\beta_1),\beta_2,\ldots,T^{s_2}(\beta_2),\ldots,\beta_t,\ldots,T^{s_t}(\beta_t)),
\end{equation}

We define $\lo= \sum_{\alpha\in \Delta} \R h_{\alpha}$, and
consider the weight spaces of $\Lambda^2(\g)$ with the order given
by $\Delta$. Let $\gamma\in \Delta$ and define
$\Gamma(\gamma)=\{\alpha \in \widehat\Gamma_1: \alpha+\gamma \in
\Phi\}$. Take $\gamma \in \Delta, \gamma \not\in \Gamma_1$ such
that $\Gamma(\gamma)\not=\emptyset$ (it is clear that such
$\gamma$ exists, otherwise $\g$ would be of type $A_1$ and in this
case $\Gamma_1 = \Gamma_2= \emptyset$). Take ${\gamma_1}$ in
$\Gamma(\gamma)$ and ${\gamma_2} \in \widehat\Gamma_2$ such that
${\gamma_1}\prec{\gamma_2}$. We will prove that the weight
${\gamma_1}+{\gamma_2}+\gamma$ occurs in $[e_\gamma - d_\gamma
e_{-\mu \gamma},\widetilde{\ro}]$, so the last one is not zero.
Now,
\begin{equation} \label{eqA}
[e_\gamma, e_{{\gamma_1}} \wedge e_{{\gamma_2}}]=
N_{\gamma,{\gamma_1}}e_{\gamma+{\gamma_1}} \wedge e_{{\gamma_2}} +
N_{\gamma,{\gamma_2}}e_{{\gamma_1}} \wedge e_{\gamma+{\gamma_2}},
\end{equation}
and, at least, the first term is not 0, because $\gamma_1º \in
\Gamma(\gamma)$. We will see that the term $e_{\gamma+{\gamma_1}}
\wedge e_{{\gamma_2}}$ is not cancelled by the other terms of
$[e_\gamma - d_\gamma e_{-\mu \gamma},\widetilde{\ro}]$, that are
of the form:
\begin{alignat}2
&[e_\gamma, e_{{\gamma'_1}} \wedge e_{{\gamma'_2}}]&=&
N_{\gamma,{\gamma'_1}}e_{\gamma+{\gamma'_1}} \wedge
e_{{\gamma'_2}} + N_{\gamma,{\gamma'_2}}e_{{\gamma'_1}} \wedge
e_{\gamma+{\gamma'_2}},  \label{eq0}\\
&[e_\gamma, e_{{-\gamma'_1}} \wedge e_{{\gamma'_2}}]&=&
N_{\gamma,{-\gamma'_1}}e_{\gamma-{\gamma'_1}} \wedge
e_{{\gamma'_2}} + N_{\gamma,{\gamma'_2}}e_{{-\gamma'_1}} \wedge
e_{\gamma+{\gamma'_2}},  \label{eq1}\\
&[e_\gamma, e_{{\gamma'_1}} \wedge e_{{-\gamma'_2}}]&=&
N_{\gamma,{\gamma'_1}}e_{\gamma+{\gamma'_1}} \wedge
e_{{-\gamma'_2}} + N_{\gamma,{-\gamma'_2}}e_{{\gamma'_1}} \wedge
e_{\gamma-{\gamma'_2}},  \label{eq2}\\
&[e_\gamma, e_{{-\gamma'_1}} \wedge e_{{-\gamma'_2}}]&=&
N_{\gamma,{-\gamma'_1}}e_{\gamma-{\gamma'_1}} \wedge
e_{{-\gamma'_2}} + N_{\gamma,{-\gamma'_2}}e_{{-\gamma'_1}} \wedge
e_{\gamma-{\gamma'_2}},  \label{eq3}\\
&[e_{-\mu(\gamma)}, e_{{\gamma'_1}} \wedge e_{{\gamma'_2}}]&=&
N_{{-\mu(\gamma)},{\gamma'_1}}e_{{-\mu(\gamma)}+{\gamma'_1}}
\wedge e_{{\gamma'_2}} +
N_{{-\mu(\gamma)},{\gamma'_2}}e_{{\gamma'_1}} \wedge
e_{{-\mu(\gamma)}+{\gamma'_2}},  \label{eq4}\\
&[e_{-\mu(\gamma)}, e_{{-\gamma'_1}} \wedge e_{{\gamma'_2}}]&=&
N_{{-\mu(\gamma)},{-\gamma'_1}}e_{{-\mu(\gamma)}-{\gamma'_1}}
\wedge e_{{\gamma'_2}} +
N_{{-\mu(\gamma)},{\gamma'_2}}e_{{-\gamma'_1}} \wedge
e_{{-\mu(\gamma)}+{\gamma'_2}},  \label{eq5}\\
&[e_{-\mu(\gamma)}, e_{{\gamma'_1}} \wedge e_{{-\gamma'_2}}]&=&
N_{{-\mu(\gamma)},{\gamma'_1}}e_{{-\mu(\gamma)}+{\gamma'_1}}
\wedge e_{{-\gamma'_2}} +
N_{{-\mu(\gamma)},{-\gamma'_2}}e_{{\gamma'_1}} \wedge
e_{{-\mu(\gamma)}-{\gamma'_2}},  \label{eq6}\\
&[e_{-\mu(\gamma)}, e_{{-\gamma'_1}} \wedge e_{{-\gamma'_2}}]&=&
N_{{-\mu(\gamma)},{-\gamma'_1}}e_{{-\mu(\gamma)}-{\gamma'_1}}
\wedge e_{{-\gamma'_2}} +
N_{{-\mu(\gamma)},{-\gamma'_2}}e_{{-\gamma'_1}} \wedge
e_{{-\mu(\gamma)}-{\gamma'_2}},  \label{eq7}
\end{alignat}

It is clear that $e_{\gamma+{\gamma_1}} \wedge e_{{\gamma_2}}$
cannot be cancelled with terms of type (\ref{eq1}), (\ref{eq2}),
(\ref{eq3}), (\ref{eq5}), (\ref{eq6}) and (\ref{eq7}) because of
there are negative roots in the factors of the terms of the RHS of
these equations (recall that $\gamma$ and $\mu(\gamma)$ are
simple).

It could be cancelation  between a term of (\ref{eq0})   and
$e_{\gamma+{\gamma_1}} \wedge e_{{\gamma_2}}$ if
\begin{enumerate}
\item[(a)] $\gamma_1 + \gamma = \gamma'_1 + \gamma$ and $\gamma_2
= \gamma'_2$, or \item[(b)] $\gamma_1 + \gamma = \gamma'_2$ and
$\gamma_2 = \gamma'_1 + \gamma$, or \item[(c)] $\gamma_1 + \gamma
= \gamma'_1$ and $\gamma_2 = \gamma'_2 + \gamma$, or \item[(d)]
$\gamma_1 + \gamma = \gamma'_2 + \gamma$ and $\gamma_2 =
\gamma'_1$.
\end{enumerate}
Now, (a) could be satisfied if and only if $\gamma_1 = \gamma'_1$
and $\gamma_2 = \gamma'_2$. In the case (b) we have that $\gamma_1
+ \gamma = \gamma'_2$, thus $\gamma \in \Gamma_2$.  As $\gamma_1
\prec \gamma_2$ and $\gamma'_1 \prec \gamma'_2$, there exist $k,s
\in \mathbb N$ such that  $T^k(\gamma_1) = \gamma_2 =
\gamma'_1+\gamma$ and $T^s(\gamma'_1) = \gamma'_2 =
\gamma_1+\gamma$. If we think in terms of the base
(\ref{ordenraices}) it is easy to see that the last two equalities
can not happen simultaneously. The case (c) results in a
contradiction because $\gamma_1 + \gamma = \gamma'_1$ implies that
$\gamma$ belongs to $\Gamma_1$. Finally, in the case (d) we have
that  $T^k(\gamma_1) = \gamma_2 = \gamma'_1$ and $T^s(\gamma'_1) =
\gamma'_2 = \gamma_1$, thus $T^{k+s}(\gamma_1)  = \gamma_1$, which
contradicts the nilpotency of $T$.

Now, it could be cancellation  between a term of (\ref{eq4}) and
$e_{\gamma+{\gamma_1}} \wedge e_{{\gamma_2}}$ if
\begin{enumerate}
\item[(a)] $\gamma_1  +  \gamma = \gamma'_1 {-\mu(\gamma)}$ and
$\gamma_2 = \gamma'_2$, or \item[(b)] $\gamma_1  +  \gamma =
\gamma'_2$ and $\gamma_2 = \gamma'_1 {-\mu(\gamma)}$, or
\item[(c)] $\gamma_1  +  \gamma = \gamma'_1$ and $\gamma_2 =
\gamma'_2 {-\mu(\gamma)}$, or \item[(d)] $\gamma_1  +  \gamma =
\gamma'_2 {-\mu(\gamma)}$ and $\gamma_2 = \gamma'_1$.
\end{enumerate}
In the case (a), as $\gamma$ and $\mu(\gamma)$ not belongs to
$\Gamma_1$, we have that $\gamma'_1 - \mu(\gamma)$ is not a root,
a contradiction because $\gamma'_1 - \mu(\gamma) = \gamma_1  +
\gamma$. In (b), as $\gamma'_1 {-\mu(\gamma)} = \gamma_2$ is a
root, we have that $\mu(\gamma)$ is in $\Gamma_1$, thus $\gamma$
is in $\Gamma_1$, a contradiction. In (c),  $\gamma_1  +  \gamma =
\gamma'_1$ implies that $\gamma \in \Gamma_1$, a contradiction. In
(d), if $\gamma \not= \mu(\gamma)$, and in consequence  $\gamma,
\mu(\gamma) \not\in \Gamma_1 \cup \Gamma_2$, we have that
$\gamma_1  +  \gamma + {\mu(\gamma)}= \gamma'_2 $, thus $\gamma
\in \Gamma_2$, a contradiction. Thus $\gamma= \mu(\gamma)$ and we
have $\gamma_1  +  \gamma = \gamma'_2 {-\gamma}$ and $\gamma_2 =
\gamma'_1$. Now, $T^k(\gamma_1) = \gamma_2 = \gamma'_1$ and
$T^s(\gamma'_1) = \gamma'_2$, so $T^{k+s}(\gamma_1) = \gamma'_2 =
\gamma_1  +  2\gamma$. If we think in terms of the base
(\ref{ordenraices}) it is easy to see that the last equality can
not happen.
\end{proof}

\medbreak

\medbreak
\subsection{Case $\sigma= \jid$.}
In this case necessarily $\Gamma_1  = \Gamma_2  = \emptyset$. See
Proposition \ref{AJancsa}.

\begin{proposition}\label{invo4} If $\ko$ is coideal of $\go$ then
$\#J = \#\Delta-1$;  say $\{\alpha\}= \Delta-J$. Assuming this,
$\ko$ is coideal of $\go$ if and only if the coefficient of
$\alpha$ in the largest root of $\Phi$ is $1$. This happens if and
only if

\begin{enumerate}
 \item $\go$ is of type $A_n$ and $\alpha\in \Delta$ is arbitrary, or

\medbreak \item $\go$ is of type $B_n$ and $\alpha$ is the
leftmost extreme
 of the Dynkin diagram ($\alpha$ is the shortest simple root), or

\medbreak \item $\go$ is of type $C_n$ and $\alpha$ is the rightmost
extreme of the Dynkin diagram ($\alpha$ is the longest simple root),
or

\medbreak \item $\go$ is of type $D_n$ and $\alpha$ is an extreme
 of the Dynkin diagram, or

\medbreak \item $\go$ is of type $E_6$ and $\alpha$ is extreme
 of the long branch of the Dynkin diagram, or

\medbreak \item $\go$ is of type $E_7$ and $\alpha$ is the extreme
 of the long branch of the Dynkin diagram.
\end{enumerate}
\end{proposition}

\begin{proof} Let us assume that $\sigma= \jid$. Recall from \eqref{jid-deltau}
that
$$\widetilde{\ro} = t \sum_{\alpha\in \Phi^+} (1 +
(-1)^{\chi_\I(\alpha)}) e_{-\alpha} \wedge e_\alpha.$$

\medbreak\begin{step}\label{pasocero}
$[h_{\gamma},\widetilde{r_0}] =0 $ for all  $\gamma \in \Delta$.
\end{step}

This is evident. Next, we compute:

\medbreak\begin{step} If $\gamma\in \Phi^+$, then

\begin{align}\label{blemmaegamma1}
\ad e_{\gamma} (\widetilde{\ro}) &= t\Big( (1 +
(-1)^{\chi_\I(\gamma)}) h_\gamma \wedge e_\gamma + \sum_{\alpha\in
\Phi^+,\gamma-\alpha \in \Phi^+} (1 + (-1)^{\chi_\I(\alpha)})
N_{\gamma,-\alpha} \, e_{\gamma-\alpha} \wedge e_\alpha\\
&\qquad+ (1 + (-1)^{\chi_\I(\gamma)}) \sum_{\alpha\in \Phi^+,
\gamma+\alpha\in \Phi^+} N_{\gamma,-\alpha}(-1)^{\chi_\I(\alpha)}
\, e_{-\alpha} \wedge e_{\gamma+\alpha}\Big),\nonumber
 \\\label{blemmaegamma2} \ad e_{-\gamma} (\widetilde{\ro})
&=t\Big( -(1 + (-1)^{\chi_\I(\gamma)}) e_{-\gamma} \wedge h_\gamma
+ \sum_{\alpha\in \Phi^+, \gamma- \alpha \in \Phi^+} (1 +
(-1)^{\chi_\I(\alpha)})
N_{-\gamma,\alpha} \, e_{-\alpha} \wedge e_{\alpha-\gamma}\\
&\qquad - (1 + (-1)^{\chi_\I(\gamma)}) \sum_{\alpha\in \Phi^+,
\alpha - \gamma\in \Phi^+} N_{-\gamma,\alpha} (-1)^{\chi_\I(\alpha
- \gamma)} \, e_{-\alpha} \wedge e_{-\gamma+\alpha}\Big).
\nonumber
\end{align}
\end{step}

\begin{proof}  Clearly,
\begin{alignat*}2
\ad e_{\gamma} (\widetilde{\ro})&&&  =t\sum_{\alpha\in
\Phi^+,\alpha\not=\gamma}
(1+(-1)^{\chi_\I(\alpha)})\big(N_{\gamma,-\alpha}e_{\gamma-\alpha}
\wedge e_\alpha +N_{\gamma,\alpha}e_{-\alpha} \wedge
e_{\gamma+\alpha}\big) + t(1+(-1)^{\chi_\I(\gamma)})h_\gamma
\wedge e_\gamma .
\end{alignat*}

Since  $\gamma \in \Phi^+$, we have
\begin{align*}
\sum_{\alpha\in \Phi^+,\alpha\not=\gamma}
&(1+(-1)^{\chi_\I(\alpha)})\big(N_{\gamma,-\alpha}e_{\gamma-\alpha}
\wedge e_\alpha +N_{\gamma,\alpha}e_{-\alpha} \wedge
e_{\gamma+\alpha}\big)\\ &= \sum_{\alpha\in \Phi^+,\gamma-\alpha
\in \Phi^+}
(1+(-1)^{\chi_\I(\alpha)})N_{\gamma,-\alpha}e_{\gamma-\alpha}
\wedge e_\alpha \\
&+ \sum_{\alpha\in \Phi^+,\gamma-\alpha \in \Phi^-}
(1+(-1)^{\chi_\I(\alpha)})N_{\gamma,-\alpha}e_{\gamma-\alpha}
\wedge
e_\alpha \\
&+\sum_{\alpha\in \Phi^+, \gamma+\alpha\in \Phi^+}
(1+(-1)^{\chi_\I(\alpha)})N_{\gamma,\alpha}e_{-\alpha} \wedge
e_{\gamma+\alpha}\\
&= \sum_{\alpha\in \Phi^+,\gamma-\alpha \in \Phi^+}
(1+(-1)^{\chi_\I(\alpha)})N_{\gamma,-\alpha}e_{\gamma-\alpha}
\wedge
e_\alpha \\
&+\sum_{\alpha\in \Phi^+, \gamma+\alpha\in \Phi^+}
\big[(1+(-1)^{\chi_\I(\alpha+\gamma)})N_{\gamma,-\alpha-\gamma} +
(1+(-1)^{\chi_\I(\alpha)}) N_{\gamma,\alpha}\big] e_{-\alpha}
\wedge e_{\gamma+\alpha}.
\end{align*}

Now,  $(-1)^{\chi_\I(\alpha+\gamma)}=-(-1)^{\chi_\I(\alpha)}
(-1)^{\chi_\I(\gamma)}$;  and $N_{\gamma,\alpha}=
N_{-\alpha-\gamma,\gamma}=-N_{\gamma,-\alpha-\gamma}$ by
\eqref{remarknalphabeta-dos} and \eqref{remarknalphabeta-tres}.
Thus the second sum in the last expression equals
$$
\sum_{\alpha\in \Phi^+, \gamma+\alpha\in \Phi^+}
N_{\gamma,\alpha}(-1)^{\chi_\I(\alpha)}(1 +
(-1)^{\chi_\I(\gamma)})e_{-\alpha} \wedge e_{\gamma+\alpha},
$$
and \eqref{blemmaegamma1} follows. The proof of
\eqref{blemmaegamma2} is completely analogous.
\end{proof}

\medbreak
\begin{step}\label{propk0}
$\ko$ is a coideal of $\go$ if and only if for any  $\gamma\in
\Phi^+$ such that $(-1)^{\chi_\I(\gamma)}= -1$,  and for any
$\alpha,\beta \in \Phi^+$ such that $\gamma=\alpha+\beta$, one has
$(-1)^{\chi_\I(\alpha)}=(-1)^{\chi_\I(\beta)}= -1$.
\end{step}

\begin{proof} By Lemma \ref{coideal} (\ref{kacaso2}) and Step \ref{pasocero},
$\ko$ is a coideal of $\go$ if and only if   $\big(1 -
(-1)^{\chi_\I(\gamma)}\big)\ad e_\gamma ( \widetilde{r_0}) = 0$
for all $\gamma\in \Phi$.   If $(-1)^{\chi_\I(\gamma)}=0$, then
there is nothing to prove. If $\gamma\in \Phi^+$ and
$(-1)^{\chi_\I(\gamma)}=1$, then by \eqref{blemmaegamma1} we have
\begin{equation*}
\big(1 - (-1)^{\chi_\I(\gamma)}\big)\ad e_{\gamma}
(\widetilde{\ro}) = 2t \sum_{\alpha\in \Phi^+,\gamma-\alpha \in
\Phi^+} (1 + (-1)^{\chi_\I(\alpha)}) N_{\gamma,-\alpha} \,
e_{\gamma-\alpha} \wedge e_\alpha.
\end{equation*}
Let $\alpha_1,\ldots,\alpha_r \in \Phi^+$ a maximal set satisfying
$\beta_1:= \gamma-\alpha_1,\ldots,\beta_r:= \gamma-\alpha_r \in
\Phi^+$ and $\alpha_i \not = \beta_j$ for all $i,j$. Then the last
equation is equivalent to:
\begin{equation*}
\big(1 - (-1)^{\chi_\I(\gamma)}\big)\ad e_{\gamma}
(\widetilde{\ro}) = 2t \sum_{i = 1}^r (1 +
(-1)^{\chi_\I(\alpha_i)}) N_{\gamma,-\alpha_i} \, e_{\beta_i}
\wedge e_{\alpha_i}+ (1 +(-1)^{\chi_\I(\beta_i)})
N_{\gamma,-\beta_i} \, e_{\alpha_i} \wedge e_{\beta_i}.
\end{equation*}
From (\ref{remarknalphabeta-dos}) and
(\ref{remarknalphabeta-tres}) we have
$$
N_{\gamma,-\alpha_i} = N_{-\beta_i,\gamma} = -N_{\gamma,-\beta_1},
$$
thus
\begin{align*}
\big(1 - (-1)^{\chi_\I(\gamma)}\big)\ad e_{\gamma}
(\widetilde{\ro}) &= 2t \sum_{i = 1}^r (1 +
(-1)^{\chi_\I(\alpha_i)}) N_{\gamma,-\alpha_i} \, e_{\beta_i}
\wedge e_{\alpha_i}\\  &\qquad\qquad+ (1 +(-1)^{\chi_\I(\beta_i)})
(-N_{\gamma,-\alpha_i}) \, e_{\alpha_i} \wedge e_{\beta_i}\\
&= 2t \sum_{i = 1}^r (2 + (-1)^{\chi_\I(\alpha_i)}+
(-1)^{\chi_\I(\beta_i)}) N_{\gamma,-\alpha_i} \, e_{\beta_i}
\wedge e_{\alpha_i}.
\end{align*}

 Then $\big(1 - (-1)^{\chi_\I(\gamma)}\big)\ad e_{\gamma}
(\widetilde{\ro}) = 0$ if and only if $(-1)^{\chi_\I(\alpha)} =
-1$ for all $\alpha\in \Phi^+$ such that $\gamma-\alpha \in
\Phi^+$. This proves the claim.
\end{proof}

\medbreak\begin{step}\label{pmenor1} If $\ko$ is a coideal of
$\go$ then $\#J = \#\Delta-1$.
\end{step}

\begin{proof} If $J = \Delta$ then $\go$ is compact, contrary to
our assumptions. Thus there is at least one element in $\Delta -
J$. Assume that there is more than one element in $\Delta - J$. We
can then choose $\alpha\neq \beta \in \Delta - J$ such that the
minimal path from $\alpha$ to $\beta$ in the Dynkin diagram
contains only points in $J$. It follows that there exists $\gamma
\in \Phi^+$ satisfying $$\gamma = \alpha + k_1\alpha_1 + \dots +
k_s\alpha_s + \beta,$$ with $\alpha_1, \dots, \alpha_s\in J$ and
$\alpha + k_1\alpha_1 + \dots + k_s\alpha_s\in \Phi^+$. Then, by
Lemma \ref{involuciones} (c):
$$
(-1)^{\chi_\I(\gamma)} = (-1)^{\chi_\J(\gamma) + \het(\gamma) + 1}
= (-1)^{2k_1  + \dots + 2k_s + 1} = -1,
$$
but $\beta \notin J$, contradicting Step \ref{propk0}.
\end{proof}

\medbreak
\begin{step}\label{cor1}
Assume that $\Delta-J = \{\alpha\}$. If $\gamma \in \Phi^+$, write
$\gamma= \sum_{\beta\in \Delta} k_\beta \beta$. Then $\ko$ is
coideal of $\go$ if and only if the coefficient $k_\alpha$ is $0$
or $1$ for any $\gamma\in\Phi^+$.
\end{step}
\begin{proof}
Assume that $\ko$ is a coideal of $\go$. If for some $\gamma\in
\Phi^+$, $k_\alpha \ge 2$ then we can assume that $k_\alpha = 2$
(for some other positive root, say). Computing
$(-1)^{\chi_\I(\gamma)}$ as in the previous step we get a
contradiction. Conversely, assume that the coefficient $k_\alpha$
is $0$ or $1$ for any $\gamma\in\Phi^+$. Note that
$(-1)^{\chi_\I(\gamma)} = -(-1)^{k_\alpha}$. Thus
$(-1)^{\chi_\I(\gamma)} = -1$ if and only if $k_\alpha = 0$. We
conclude now from Step \ref{propk0}.
\end{proof}

 Recall that the {\it largest root} of $\Phi$ is the
highest weight of the adjoint representation of $\g$.

\medbreak\begin{step}\label{cor2} Assume that $\Delta-J =
\{\alpha\}$.
 Then $\ko$ is coideal of $\go$ if and
only if the coefficient of $\alpha$ in the largest root of $\Phi$
is $1$.
\end{step}

\begin{proof}
If $\sum_{\beta \in \Delta} t_\beta \beta$ is the largest root and
$\gamma = \sum_{\beta \in \Delta} k_\beta \beta\in \Phi^+$, then
 $k_\beta \le t_\beta$. Thus the Step follows
immediately from Step \ref{cor1}.
\end{proof}

It remains only to determine the Dynkin diagrams with a simple
root whose coefficient in the largest root is $1$. This is an easy
task, by inspecting the largest root of each system as listed in
\cite{Kn}, Appendix C. For example, the largest root corresponding
to a system of type $B_n$ is $\alpha_1 + \sum_{i=2}^n 2\alpha_i$.
So, from Step \ref{cor2}, $\ko$ is subcoideal of $\go$ if and only
if $\alpha=\alpha_1$. The same argument applies to the other
Dynkin diagrams.
\end{proof}

\medbreak
\subsection{Case $\sigma= \jdynkin$, $\mu\not=\id$,
$\Gamma_1  = \Gamma_2  = \emptyset$.}

We shall show that $\ko$ is not coideal of $\go$ except for $\go$
of type $A_2$. We begin by the following reduction; recall the set
$\I$ defined in Lemma \ref{involuciones}.

\begin{lemma}\label{ugamma} If $\gamma \in \Phi^+$ then
$\ad(e_\gamma  - (-1)^{\chi_\I(\gamma)} e_{\mu(\gamma)})
(\widetilde{\ro}) \equiv t u_\gamma \mod \h \otimes \g + \g
\otimes \h$, where
\begin{align*}
u_\gamma = &\quad\sum_{\alpha\in \Phi^+, \gamma - \alpha\in
\Phi^+}
N_{\gamma,-\alpha}e_{\gamma-\alpha} \wedge e_\alpha \\
 &+ \sum_{\alpha\in \Phi^+, \gamma - \alpha\in \Phi^+}
(-1)^{\chi_\I(\gamma)+1}
N_{\mu(\gamma),-\mu(\alpha)}e_{\mu(\gamma-\alpha)} \wedge
e_{\mu(\alpha)} \\
&+\sum_{\alpha\in \Phi^+, \gamma - \alpha \in \Phi^+}
(-1)^{\chi_\I(\alpha)} (N_{\gamma,-\alpha}
-\overline{N_{-\gamma,\alpha}}) e_{\gamma-\alpha}\wedge
e_{\mu(\alpha)}\\
 &+\sum_{\alpha\in \Phi^+,  \gamma+\alpha \in \Phi^+}
(-1)^{\chi_\I(\gamma+\alpha)}( N_{\alpha,\gamma} +
\overline{N_{-\alpha,-\gamma}}) e_{-\alpha}\wedge
e_{\mu(\gamma+\alpha)}\\
&+ \sum_{\alpha\in \Phi^+,\gamma+\alpha\in \Phi^+}
(-1)^{\chi_\I(\alpha)}(N_{\gamma,\alpha}+
\overline{N_{-\gamma,-\alpha}} ) e_{-\mu(\alpha)}\wedge e_{\gamma
+\alpha}
\end{align*}
\end{lemma}

\begin{proof} For shortness, let $``X\equiv Y"$ mean
$``X \equiv Y\mod \h \otimes \g + \g \otimes \h"$. We compute:

\begin{align*}
\frac1t\ad e_{\gamma} (\widetilde{\ro})  &\equiv\, \sum_{\alpha\in
\Phi^+} \Big([e_\gamma,e_{-\alpha}] \wedge e_\alpha+ e_{-\alpha}
\wedge [e_\gamma,e_\alpha] + (-1)^{\chi_\I(\alpha)}
\big([e_\gamma,e_{-\mu (\alpha)}]\wedge e_\alpha +
e_{-\mu(\alpha)}\wedge [e_\gamma,e_\alpha]\big)\Big) \\
&\equiv\,  \sum_{\alpha\in \Phi^+, \alpha\not=\gamma}
\big(N_{\gamma,-\alpha}e_{\gamma-\alpha} \wedge e_\alpha
+ N_{\gamma,\alpha} e_{-\alpha} \wedge e_{\gamma +\alpha}\big) \\
&\qquad +  \sum_{\alpha\in \Phi^+, \alpha\not=\mu(\gamma)}
(-1)^{\chi_\I(\alpha)} N_{\gamma,-\mu(\alpha)} e_{\gamma-\mu
(\alpha)}\wedge e_\alpha  \\
&\qquad +   \,   \sum_{\alpha\in \Phi^+}
(-1)^{\chi_\I(\alpha)}N_{\gamma,\alpha} e_{-\mu(\alpha)}\wedge
e_{\gamma +\alpha}\end{align*}
\begin{align*} &\equiv\,
\\&\quad\qquad \sum_{\alpha\in \Phi^+, \gamma - \alpha\in \Phi^+}
N_{\gamma,-\alpha}e_{\gamma-\alpha} \wedge e_\alpha  \tag A\\
&\qquad + \,  \,  \sum_{\alpha\in \Phi^+, \alpha - \mu(\gamma) \in
\Phi^+} (-1)^{\chi_\I(\alpha)} N_{\gamma,-\mu(\alpha)}
e_{\gamma-\mu
(\alpha)}\wedge e_\alpha \tag B\\
&\qquad + \,  \,  \sum_{\alpha\in \Phi^+, \mu(\gamma) - \alpha \in
\Phi^+} (-1)^{\chi_\I(\alpha)} N_{\gamma,-\mu(\alpha)}
e_{\gamma-\mu (\alpha)}\wedge e_\alpha \tag C\\
&\qquad  + \,  \,  \sum_{\alpha\in \Phi^+,\gamma+\alpha\in \Phi^+}
(-1)^{\chi_\I(\alpha)}N_{\gamma,\alpha} e_{-\mu(\alpha)}\wedge
e_{\gamma +\alpha}.\tag D
\end{align*}

Here in the third congruence we use that
\begin{alignat*}2
\sum_{\alpha\in \Phi^+,  \alpha - \gamma\in \Phi^+}
N_{\gamma,-\alpha}e_{\gamma-\alpha} \wedge e_\alpha +
\sum_{\alpha\in \Phi^+, \alpha\neq\gamma}N_{\gamma,\alpha}
e_{-\alpha} \wedge e_{\gamma +\alpha} =0
\end{alignat*}
by \eqref{remarknalphabeta-dos} and \eqref{remarknalphabeta-tres}.
Changing $\alpha$ by $\mu(\alpha)$, we have
\begin{align*}
&\qquad  \,  \,  \sum_{\alpha\in \Phi^+, \alpha - \gamma \in
\Phi^+} (-1)^{\chi_\I(\mu(\alpha))} N_{\gamma,-\alpha} e_{\gamma-\alpha}\wedge e_{\mu(\alpha)} \tag B\\
&\qquad  \,  \,  \sum_{\alpha\in \Phi^+, \gamma - \alpha \in
\Phi^+} (-1)^{\chi_\I(\mu(\alpha))} N_{\gamma,-\alpha}
e_{\gamma-\alpha}\wedge e_{\mu(\alpha)}. \tag C
\end{align*}
Now,
\begin{align*}
  \sum_{\alpha\in \Phi^+, \alpha - \gamma \in
\Phi^+} &(-1)^{\chi_\I(\mu(\alpha))} N_{\gamma,-\alpha}
e_{\gamma-\alpha}\wedge e_{\mu(\alpha)} \tag B\\
&=  \sum_{\alpha\in \Phi^+,  \gamma+\alpha \in \Phi^+}
(-1)^{\chi_\I(\mu(\gamma+\alpha))} N_{\gamma,-\gamma-\alpha}
e_{-\alpha}\wedge e_{\mu(\gamma+\alpha)} \\
&=  \sum_{\alpha\in \Phi^+,  \gamma+\alpha \in \Phi^+}
(-1)^{\chi_\I(\mu(\gamma+\alpha))} N_{\alpha,\gamma}
e_{-\alpha}\wedge e_{\mu(\gamma+\alpha)}.
\end{align*}
The first equality follows from the change of variables $\alpha$
by $\alpha-\gamma$, and the second from
\eqref{remarknalphabeta-tres}. In analogous way, we have
\begin{align*}
\frac1t(-1)^{\chi_\I(\gamma)+1}\ad e_{\mu(\gamma)} (\widetilde{\ro})  \equiv \,&  \\
& \,  \sum_{\alpha\in \Phi^+, \mu(\gamma) - \alpha\in \Phi^+}
(-1)^{\chi_\I(\gamma)+1}
N_{\mu(\gamma),-\alpha}e_{\mu(\gamma)-\alpha} \wedge e_\alpha \tag E\\
\qquad +\, & \, \sum_{\alpha\in \Phi^+, \alpha - \gamma \in
\Phi^+}
 (-1)^{\chi_\I(\alpha)+\chi_\I(\gamma)+1}
N_{\mu(\gamma),-\mu(\alpha)}
e_{\mu(\gamma)-\mu(\alpha)}\wedge e_\alpha \tag F\\
\qquad + \, & \,  \sum_{\alpha\in \Phi^+, \gamma - \alpha \in
\Phi^+} (-1)^{\chi_\I(\alpha)+\chi_\I(\gamma)+1} N_{\mu(\gamma),
-\mu(\alpha)}
e_{\mu(\gamma)-\mu (\alpha)}\wedge e_\alpha \tag G\\
\qquad + \,&\, \sum_{\alpha\in \Phi^+, \mu(\gamma)+\alpha\in
\Phi^+}
(-1)^{\chi_\I(\alpha)+\chi_\I(\gamma)+1}N_{\mu(\gamma),\alpha}
e_{-\mu(\alpha)}\wedge e_{\mu(\gamma) +\alpha}. \tag H
\end{align*}

By \eqref{mualpha}, we have
\begin{equation}
(-1)^{\chi_\I(\alpha)+\chi_\I(\gamma)+1}
N_{\mu(\gamma),-\mu(\alpha)} =
(-1)^{\chi_\I(-\mu(\gamma))+\chi_\I(\mu(\alpha))+1}
N_{\mu(\gamma),-\mu(\alpha)} =
(-1)^{\chi_\I(\mu(\gamma-\alpha))+1}
\overline{N_{-\gamma,\alpha}}.
\end{equation}
Changing $\gamma - \alpha$ by $\alpha$ and using $x \wedge y = - y
\wedge x$, we have
\begin{align*}
& \,  \sum_{\alpha\in \Phi^+, \gamma - \alpha\in \Phi^+}
(-1)^{\chi_\I(\gamma)+1}
N_{\mu(\gamma),-\mu(\alpha)}e_{\mu(\gamma-\alpha)} \wedge e_{\mu(\alpha)} \tag E\\
\qquad \, & \, \sum_{\alpha\in \Phi^-, \gamma-\alpha \in \Phi^+}
  (-1)^{\chi_\I(\mu(\alpha))}
\overline{N_{-\gamma,\gamma-\alpha}}
e_{\gamma-\alpha} \wedge e_{\mu(\alpha)}\tag F\\
\qquad  \, & \,  \sum_{\alpha\in \Phi^+, \gamma - \alpha \in
\Phi^+} (-1)^{\chi_\I(\mu(\alpha))}
\overline{N_{-\gamma,\gamma-\alpha}} e_{\gamma-\alpha} \wedge
e_{\mu (\alpha)}\tag G\end{align*}

Now, performing the change of variable $\alpha$ by $-\alpha$ and
using \eqref{remarknalphabeta-dos} and
\eqref{remarknalphabeta-tres}, we have
\begin{align*}
\qquad \, & \, \sum_{\alpha\in \Phi^+, \gamma+\alpha \in \Phi^+}
  (-1)^{\chi_\I(\mu(\alpha))}
\overline{N_{-\gamma,\gamma+\alpha}} e_{\gamma+\alpha} \wedge
e_{-\mu(\alpha)} \tag F
\\\qquad \, &= \, \sum_{\alpha\in \Phi^+,
\gamma+\alpha \in \Phi^+}
  (-1)^{\chi_\I(\mu(\alpha))}
\overline{N_{-\alpha,-\gamma}} e_{\gamma+\alpha} \wedge
e_{-\mu(\alpha)}
\\\qquad \, &= \, \sum_{\alpha\in \Phi^+, \gamma+\alpha \in \Phi^+}
  (-1)^{\chi_\I(\mu(\alpha))}
\overline{N_{-\gamma,-\alpha}}  e_{-\mu(\alpha)} \wedge
e_{\gamma+\alpha}.
\end{align*}

By \eqref{remarknalphabeta-tres}, we have the following expression
for (G):
\begin{align*}
\sum_{\alpha\in \Phi^+, \gamma - \alpha \in \Phi^+}
(-1)^{\chi_\I(\mu(\alpha))+1} \overline{N_{-\gamma,\alpha}}
e_{\gamma-\alpha} \wedge e_{\mu (\alpha)}\tag G
\end{align*}

For (H), we perform the change of variables $\alpha$ by
$\mu(\alpha)$; applying \eqref{chimu},
\eqref{remarknalphabeta-dos}, we get:
\begin{align*}
\qquad  \sum_{\alpha\in \Phi^+, \gamma+\alpha\in \Phi^+}&
(-1)^{\chi_\I(\mu(\alpha))+\chi_\I(\gamma)+1}N_{\mu(\gamma),\mu(\alpha)}
e_{-\alpha}\wedge e_{\mu(\gamma+\alpha)} \tag H \\
&= \sum_{\alpha\in \Phi^+, \gamma+\alpha\in \Phi^+}
(-1)^{\chi_\I(-\mu(\gamma))+\chi_\I(-\mu(\alpha))+1}N_{\mu(\gamma),\mu(\alpha)}
e_{-\alpha}\wedge e_{\mu(\gamma+\alpha)} \\
&= \sum_{\alpha\in \Phi^+, \gamma+\alpha\in \Phi^+}
(-1)^{\chi_\I(\mu(\gamma+\alpha))}\overline{N_{-\alpha,-\gamma}}
e_{-\alpha}\wedge e_{\mu(\gamma+\alpha)}.
\end{align*}

Finally, (B) $+$ (H) is
\begin{equation*}
 \sum_{\alpha\in \Phi^+,  \gamma+\alpha \in \Phi^+}
(-1)^{\chi_\I(\mu(\gamma+\alpha))}( N_{\alpha,\gamma} +
\overline{N_{-\alpha,-\gamma}}) e_{-\alpha}\wedge
e_{\mu(\gamma+\alpha)},
\end{equation*}
(C) $+$ (G) is
\begin{equation*}
\sum_{\alpha\in \Phi^+, \gamma - \alpha \in \Phi^+}
(-1)^{\chi_\I(\mu(\alpha))} (N_{\gamma,-\alpha}
-\overline{N_{-\gamma,\alpha}}) e_{\gamma-\alpha}\wedge
e_{\mu(\alpha)},
\end{equation*}
and (D) $+$ (F) is
\begin{equation*}
\sum_{\alpha\in \Phi^+,\gamma+\alpha\in \Phi^+}
(-1)^{\chi_\I(\alpha)}(N_{\gamma,\alpha}+
\overline{N_{-\gamma,-\alpha}} ) e_{-\mu(\alpha)}\wedge e_{\gamma
+\alpha}
\end{equation*}
\end{proof}

\begin{remark}
If $\gamma \in \Phi^+$ then by the previous Lemma,
$$\ad(e_\gamma - (-1)^{\chi_\I(\gamma)} e_{\mu(\gamma)})
(\widetilde{\ro}) \equiv t v_\gamma \mod \big(\h \otimes \g + \g
\otimes \h + \n^+ \otimes \n^- + \n^- \otimes \n^+\big),$$ where
\begin{align*}
v_\gamma = &\quad\sum_{\alpha\in \Phi^+, \gamma - \alpha\in
\Phi^+}
N_{\gamma,-\alpha}e_{\gamma-\alpha} \wedge e_\alpha \\
 &+ \sum_{\alpha\in \Phi^+, \gamma - \alpha\in \Phi^+}
(-1)^{\chi_\I(\gamma)+1}
N_{\mu(\gamma),-\mu(\alpha)}e_{\mu(\gamma-\alpha)} \wedge
e_{\mu(\alpha)} \\
&+\sum_{\alpha\in \Phi^+, \gamma - \alpha \in \Phi^+}
(-1)^{\chi_\I(\alpha)} (N_{\gamma,-\alpha}
-\overline{N_{-\gamma,\alpha}}) e_{\gamma-\alpha}\wedge
e_{\mu(\alpha)}.
\end{align*}
It is clear that $v_\gamma \not=0$ implies $\ad(e_\gamma -
(-1)^{\chi_{\I}(\gamma)} e_{\mu(\gamma)})
(\widetilde{\ro})\not=0$.
\end{remark}

Let $V$ a finite dimensional representation of $\g$. If $\lambda
\in \h^*$, then we denote $V_{(\lambda)}:= \{v \in V: h.v
=\lambda(h)v\}$. If $v \in V$, then we denote $v_{(\lambda)}$ the
component of weight $\lambda$ of $v$.

\begin{corollary} \label{cor:1.9}
If $\gamma \in \Phi^+$ then
$$\ad(e_\gamma - (-1)^{\chi_\I(\gamma)} e_{\mu(\gamma)})
(\widetilde{r_0})_{(\gamma)} \equiv t w_\gamma \mod \big(\h
\otimes \g + \g \otimes \h + \n^+ \otimes \n^- + \n^- \otimes
\n^+\big),$$ where
\begin{align*}
w_\gamma = \, 2\sum_{\alpha\in \Phi^+, \gamma - \alpha\in \Phi^+}
N_{\gamma,-\alpha}e_{\gamma-\alpha} \wedge e_\alpha
+\sum_{\alpha\in (\Phi^+)^\mu, \gamma - \alpha \in \Phi^+}
(-1)^{\chi_\I(\alpha)} (N_{\gamma,-\alpha}
-\overline{N_{-\gamma,\alpha}}) e_{\gamma-\alpha}\wedge
e_{\alpha}.
\end{align*}\qed
\end{corollary}

\begin{lemma}\label{lemma:a3}
Let $\sigma= \jdynkin$, with $J= \{\alpha_2\}$. Let $\alpha_1,
\alpha_3 \in \Delta$ such that $\alpha_1,\alpha_2,\alpha_3$ is a
subsystem of type $A_3$ and $\sigma(\alpha_1)=\alpha_3$. Then
$\chi_\I(\alpha_1+\alpha_2+\alpha_3)=1$.
\end{lemma}
\begin{proof}
By the Jacobi identity we have that
$$
[e_{\alpha_1}, [e_{\alpha_2},e_{\alpha_3}]] = - [e_{\alpha_2},
[e_{\alpha_3},e_{\alpha_1}]]-[e_{\alpha_3},
[e_{\alpha_1},e_{\alpha_2}]]  =-[e_{\alpha_3},
[e_{\alpha_1},e_{\alpha_2}]].
$$
Now $[e_{\alpha_1}, [e_{\alpha_2},e_{\alpha_3}]] =
N_{\alpha_2,\alpha_3} N_{\alpha_1,\alpha_2+\alpha_3}$ and
$-[e_{\alpha_3}, [e_{\alpha_1},e_{\alpha_2}]]=
-N_{\alpha_1,\alpha_2} N_{\alpha_3,\alpha_1+\alpha_1} $, thus
\begin{equation}\label{eq:nab1}
N_{\alpha_2,\alpha_3} N_{\alpha_1,\alpha_2+\alpha_3} =
-N_{\alpha_1,\alpha_2} N_{\alpha_3,\alpha_1+\alpha_1}.
\end{equation}

From \eqref{mualpha}, we have that
$$
(-1)^{\chi_\I(\alpha_2)+\chi_\I(\alpha_3)}N_{\alpha_2, \alpha_1} =
(-1)^{\chi_\I(\alpha_2 +\alpha_3)}
\overline{N_{-\alpha_2,-\alpha_3}}.
$$
As $\chi_\I(\alpha_2)=1$ and $\chi_\I(\alpha_3)=0$, we have
\begin{equation}\label{eq:nab2}
N_{\alpha_1, \alpha_2}= -N_{\alpha_2, \alpha_1} =
(-1)^{\chi_\I(\alpha_2 +\alpha_3)}
\overline{N_{-\alpha_2,-\alpha_3}}.
\end{equation}

Again from \eqref{mualpha}, we have
$$
(-1)^{\chi_\I(\alpha_1)+\chi_\I(\alpha_2+\alpha_3)}N_{\alpha_3,
\alpha_2+\alpha_1} = (-1)^{\chi_\I(\alpha_1+\alpha_2 +\alpha_3)}
\overline{N_{-\alpha_1,-\alpha_2-\alpha_3}}.
$$
As  $\chi_\I(\alpha_1)=0$, we get
\begin{equation}\label{eq:nab3}
(-1)^{\chi_\I(\alpha_2+\alpha_3)}N_{\alpha_3, \alpha_2+\alpha_1} =
(-1)^{\chi_\I(\alpha_1+\alpha_2 +\alpha_3)}
\overline{N_{-\alpha_1,-\alpha_2-\alpha_3}}.
\end{equation}

We  multiply equation (\ref{eq:nab3}) by $N_{\alpha_1,\alpha_2}$
and obtain
$$
(-1)^{\chi_\I(\alpha_2+\alpha_3)}N_{\alpha_1,\alpha_2}N_{\alpha_3,
\alpha_2+\alpha_1} = (-1)^{\chi_\I(\alpha_1+\alpha_2 +\alpha_3)}
N_{\alpha_1,\alpha_2}\overline{N_{-\alpha_1,-\alpha_2-\alpha_3}},
$$
applying \eqref{eq:nab1} to the $N_{\alpha_1,\alpha_2}$ on the
left, we have
$$
-(-1)^{\chi_\I(\alpha_2+\alpha_3)}N_{\alpha_2,\alpha_3}
N_{\alpha_1,\alpha_2+\alpha_3} = (-1)^{\chi_\I(\alpha_1+\alpha_2
+\alpha_3)}
N_{\alpha_1,\alpha_2}\overline{N_{-\alpha_1,-\alpha_2-\alpha_3}}.
$$
Applying (\ref{eq:nab2}) we obtain
$$
-(-1)^{\chi_\I(\alpha_2+\alpha_3)}N_{\alpha_2,\alpha_3}
N_{\alpha_1,\alpha_2+\alpha_3} = (-1)^{\chi_\I(\alpha_1+\alpha_2
+\alpha_3)} (-1)^{\chi_\I(\alpha_2 +\alpha_3)}
\overline{N_{-\alpha_2,-\alpha_3}}\overline{N_{-\alpha_1,-\alpha_2-\alpha_3}},
$$
thus
\begin{equation*} \label{eq:nab4}
-N_{\alpha_2,\alpha_3} N_{\alpha_1,\alpha_2+\alpha_3} =
(-1)^{\chi_\I(\alpha_1+\alpha_2 +\alpha_3)}
\overline{N_{-\alpha_2,-\alpha_3}}\overline{N_{-\alpha_1,-\alpha_2-\alpha_3}}.
\end{equation*}

Now, from \eqref{mualpha5} we have that $N_{-\alpha_2,-\alpha_3} =
c_1 N_{\alpha_2,\alpha_3}$ and
$N_{-\alpha_1,-\alpha_2-\alpha_3}=c_2
N_{\alpha_1,\alpha_2+\alpha_3}$ with $c_1,c_2<0$, thus
$$
- c N_{\alpha_2,\alpha_3} N_{\alpha_1,\alpha_2+\alpha_3} =
(-1)^{\chi_\I(\alpha_1+\alpha_2 +\alpha_3)}
\overline{N_{\alpha_2,\alpha_3}}\overline{N_{\alpha_1,\alpha_2+\alpha_3}}
$$
with $c>0$. This clearly implies that $\chi_\I(\alpha_1+\alpha_2
+\alpha_3)=1$.
\end{proof}

\begin{proposition}\label{invo5}
Let $\sigma= \jdynkin$, $\Gamma_1 = \Gamma_2 = \emptyset$. If $\g$
is of type $A_2$, then $\ko$ is coideal of $\go$. In the other
cases, $\ko$ is not coideal of $\go$.
\end{proposition}
\begin{proof}  When $\g$ is of type $A_2$ is easy to check directly that
$\ad(e_\gamma  - (-1)^{\chi_\I(\gamma)} e_{\mu(\gamma)})
(\widetilde{r_0})=0$ for all $\gamma \in \Phi^+$.

In the other cases, the Dynkin diagrams that admit non trivial
automorphism have subdiagrams of type $A_3$ or $A_4$ where $\mu$
acts non trivially. In the follows we prove that if we have
subdiagrams of type $A_3$ or $A_4$ where $\mu$ acts non trivially,
then there exist $\gamma \in \Phi^+$ such that $v_\gamma \not= 0$.

\medbreak

{\bf Type $A_3$}.  Let $\alpha_1$, $\alpha_2$, $\alpha_3$ be
simple roots such that they determine a subdiagram of type $A_3$
and $\mu$ restricted to $\alpha_1,\alpha_2, \alpha_3$ is non
trivial, i.e. $\mu(\alpha_1) = \alpha_3$ and  $\mu(\alpha_2) =
\alpha_2$. Now we will consider two subcases $J=\emptyset$ and
$J\not=\emptyset$.

\medbreak

{\bf $A_3$ and $J=\emptyset$}. Let $\gamma = \alpha_1+\alpha_2$,
thus $\mu(\gamma) =\alpha_2+\alpha_3$. Then
\begin{align*}
v_\gamma = &\quad
N_{\gamma,-\alpha_1}e_{\alpha_2} \wedge e_{\alpha_1} +N_{\gamma,-\alpha_2}e_{\alpha_1} \wedge e_{\alpha_2}\\
 &+ (-1)^{\chi_\I(\gamma)+1}
N_{\mu(\gamma),-\alpha_3}e_{\alpha_2} \wedge e_{\alpha_3}
+(-1)^{\chi_\I(\gamma)+1} N_{\mu(\gamma),-\alpha_2}e_{\alpha_3}
\wedge e_{\alpha_2} \\
&+ (-1)^{\chi_\I(\alpha_1)} (N_{\gamma,-\alpha_1}
-\overline{N_{-\gamma,\alpha_1}}) e_{\alpha_2}\wedge
e_{\alpha_3}+(-1)^{\chi_\I(\alpha_2)} (N_{\gamma,-\alpha_2}
-\overline{N_{-\gamma,\alpha_2}}) e_{\alpha_1}\wedge e_{\alpha_2}.
\end{align*}
Now, $N_{\gamma,-\alpha_1}= N_{\alpha_1+\alpha_2,-\alpha_1} =
N_{-\alpha_1, -\alpha_2} = N_{-\alpha_2,\alpha_1+\alpha_2 } =
-N_{\gamma,-\alpha_2}$. As $J=\emptyset$, we have
$\chi_\I(\alpha_1)=\chi_\I(\alpha_2)=0$. Finally, by
\eqref{mualpha5}, we have
$$
\overline{N_{-\gamma,\alpha_2}} =
\frac{1}{(-\gamma|\alpha_2)\overline{N_{\gamma,-\alpha_2}}}.
$$
Thus,
\begin{align*}
v_\gamma = &\quad
(3N_{\gamma,-\alpha_2}-\frac{1}{(-\gamma|\alpha_2)\overline{N_{\gamma,-\alpha_2}}})e_{\alpha_1}\wedge
e_{\alpha_2} + c e_{\alpha_2} \wedge e_{\alpha_3} \\
&=\frac{3N_{\gamma,-\alpha_2}\overline{N_{\gamma,-\alpha_2}}(\gamma|\alpha_2)+
1}{(\gamma|\alpha_2)\overline{N_{\gamma,-\alpha_2}}}e_{\alpha_1}\wedge
e_{\alpha_2} + c
e_{\alpha_2} \wedge e_{\alpha_3}  \\
&=\frac{3||N_{\gamma,-\alpha_2}||^2(\gamma|\alpha_2)+
1}{(\gamma|\alpha_2)\overline{N_{\gamma,-\alpha_2}}}e_{\alpha_1}\wedge
e_{\alpha_2} + c e_{\alpha_2} \wedge e_{\alpha_3}.
\end{align*}
As $(\gamma|\alpha_2)>0$, we have $v_\gamma\not = 0$.

\medbreak {\bf Type $A_3$ and $J\not=\emptyset$}. Let $\gamma =
\alpha_1+\alpha_2+\alpha_3$, thus $\mu(\gamma) =\gamma$ and from
Lemma \ref{lemma:a3} we have $\chi_\I(\gamma)=1$. From
Corollary 2.10 we have that $v_\gamma = w_\gamma$ andFrom Corollary
\ref{cor:1.9} we have
\begin{align*}
w_\gamma = &\quad2\sum_{\alpha\in \Phi^+, \gamma - \alpha\in
\Phi^+}
N_{\gamma,-\alpha}e_{\gamma-\alpha} \wedge e_\alpha  \\
&+\sum_{\alpha\in (\Phi^+)^\mu, \gamma - \alpha \in \Phi^+}
(-1)^{\chi_\I(\alpha)} (N_{\gamma,-\alpha}
-\overline{N_{-\gamma,\alpha}}) e_{\gamma-\alpha}\wedge e_{\alpha}
\end{align*}
Now there  exists no $\alpha\in (\Phi^+)^\mu$ such that $\gamma -
\alpha \in \Phi^+$; hence we have
\begin{align*}
\frac12w_\gamma = &\sum_{\alpha\in \Phi^+, \gamma - \alpha\in
\Phi^+,}
N_{\gamma,-\alpha}e_{\gamma-\alpha} \wedge e_\alpha  \\
= & N_{\gamma,-\alpha_1}e_{\alpha_2+\alpha_3} \wedge e_{\alpha_1}
+ N_{\gamma,-\alpha_3}e_{\alpha_1+\alpha_2} \wedge e_{\alpha_3} +
N_{\gamma,-\alpha_1-\alpha_2}e_{\alpha_3} \wedge
e_{\alpha_1+\alpha_2}+ N_{\gamma,-\alpha_2-\alpha_3}e_{\alpha_1}
\wedge e_{\alpha_2+\alpha_3} \\
= &  (N_{\gamma,-\alpha_2-\alpha_3} -N_{\gamma,-\alpha_1})
e_{\alpha_1} \wedge e_{\alpha_2+\alpha_3}
+(N_{\gamma,-\alpha_1-\alpha_2}- N_{\gamma,-\alpha_3})e_{\alpha_3}
\wedge e_{\alpha_1+\alpha_2}.
\end{align*}
Now $N_{\gamma,-\alpha_2-\alpha_3} = N_{-\alpha_1,\gamma} =
-N_{\gamma,-\alpha_1}$ and $N_{\gamma,-\alpha_1-\alpha_2} =
N_{-\alpha_3,\gamma} = -N_{\gamma,-\alpha_3}$, thus
\begin{align*}
\frac12w_\gamma =   2N_{\gamma,-\alpha_2-\alpha_3} e_{\alpha_1}
\wedge e_{\alpha_2+\alpha_3} +
2N_{\gamma,-\alpha_1-\alpha_2}e_{\alpha_3} \wedge
e_{\alpha_1+\alpha_2}  \not= 0.
\end{align*}

\bigbreak {\bf Type $A_4$}. Let $\alpha_1$, $\alpha_2$,
$\alpha_3$, $\alpha_4$ be simple roots such that they determine a
subdiagram of type $A_4$ where $\mu$ restricted to
$\alpha_1,\alpha_2, \alpha_3, \alpha_4$ is non trivial, i. e.
$\mu(\alpha_1) = \alpha_4$ and $\mu(\alpha_2) = \alpha_3$. Let
$\gamma = \alpha_1+\alpha_2$ and it is clear that $\mu(\gamma)
=\alpha_3+\alpha_4$. Then
\begin{align*}
v_\gamma = &
N_{\gamma,-\alpha_1}e_{\alpha_2} \wedge e_{\alpha_1} +N_{\gamma,-\alpha_2}e_{\alpha_1} \wedge e_{\alpha_2}\\
 &+ (-1)^{\chi_\I(\gamma)+1}
N_{\mu(\gamma),-\alpha_4}e_{\alpha_3} \wedge e_{\alpha_4} +
(-1)^{\chi_\I(\gamma)+1} N_{\mu(\gamma),-\alpha_3}e_{\alpha_4}
\wedge e_{\alpha_3}\\
&+ (-1)^{\chi_\I(\alpha_1)} (N_{\gamma,-\alpha_1}
-\overline{N_{-\gamma,\alpha_1}}) e_{\alpha_2}\wedge e_{\alpha_4}+
(-1)^{\chi_\I(\alpha_2)} (N_{\gamma,-\alpha_2}
-\overline{N_{-\gamma,\alpha_2}}) e_{\alpha_1}\wedge e_{\alpha_3}.
\end{align*}

Now, from \eqref{remarknalphabeta-dos} and
\eqref{remarknalphabeta-tres}, we have

$$
N_{\gamma,-\alpha_1}e_{\alpha_2} \wedge e_{\alpha_1}
+N_{\gamma,-\alpha_2}e_{\alpha_1} \wedge e_{\alpha_2} =
(N_{\gamma,-\alpha_2}- N_{\gamma,-\alpha_1})e_{\alpha_1} \wedge
e_{\alpha_2} = 2N_{\gamma,-\alpha_2}e_{\alpha_1} \wedge
e_{\alpha_2} \not=0.
$$
\end{proof}

\end{document}